\documentclass[11pt]{amsart}
\usepackage{amsmath,amssymb,mathrsfs}
\newtheorem{theorem}{Theorem}[section]
\newtheorem{proposition}[theorem]{Proposition}
\newtheorem{corollary}[theorem]{Corollary}
\newtheorem{lemma}[theorem]{Lemma}

\theoremstyle{definition}
\newtheorem{definition}[theorem]{Definition}
\newtheorem{example}[theorem]{Example}
\begin{document}

\title{Singularity models in the three-dimensional Ricci flow}
\author{Simon Brendle}
\address{Department of Mathematics \\ Columbia University \\ New York NY 10027}
\begin{abstract}
The Ricci flow is a natural evolution equation for Riemannian metrics on a given manifold. The main goal is to understand singularity formation. In his spectacular 2002 breakthrough, Perelman achieved a qualitative understanding of singularity formation in dimension $3$. More precisely, Perelman showed that every finite-time singularity to the Ricci flow in dimension $3$ is modeled on an ancient $\kappa$-solution. Moreover, Perelman proved a structure theorem for ancient $\kappa$-solutions in dimension $3$. 

In this survey, we discuss recent developments which have led to a complete classification of all the singularity models in dimension $3$. Moreover, we give an alternative proof of the classification of noncollapsed steady gradient Ricci solitons in dimension $3$ (originally proved by the author in 2012).
\end{abstract}
\thanks{The author was supported by the National Science Foundation under grants DMS-1806190 and DMS-2103573 and by the Simons Foundation.}

\maketitle 

\section{Background on the Ricci flow}

\label{background}

Geometric evolution equations play an important role in differential geometry. The most important such evoluation equation is the Ricci flow for Riemannian metrics which was introduced by Hamilton \cite{Hamilton1}:

\begin{definition}[R.~Hamilton \cite{Hamilton1}]
Let $g(t)$ be a one-parameter family of Riemannian metrics on a manifold $M$. We say that the metrics $g(t)$ evolve by the Ricci flow if 
\begin{equation} 
\frac{\partial}{\partial t} g(t) = -2 \, \text{\rm Ric}_{g(t)}. 
\end{equation} 
\end{definition}

In his paper \cite{Hamilton1}, Hamilton established short time existence and uniqueness for the Ricci flow. 

\begin{theorem}[R.~Hamilton \cite{Hamilton1}; D.~DeTurck \cite{DeTurck}]
\label{short.time.existence}
Let $g_0$ be a Riemannian metric on a compact manifold $M$. Then there exists a unique solution $g(t)$, $t \in [0,T)$, to the Ricci flow with initial metric $g(0)=g_0$. Here, $T$ is a positive real number which depends on the initial data.
\end{theorem}

The main difficulty in proving Theorem \ref{short.time.existence} is that the Ricci flow is weakly, but not strictly, parabolic. This is due to the fact that the Ricci flow is invariant under the diffeomorphism group of $M$. This problem can be overcome using DeTurck's trick \cite{DeTurck}. In the following, we sketch the argument (see \cite{Brendle2} or \cite{Topping} for details). Let us fix a compact manifold $M$ and a smooth one-parameter family of background metrics $h(t)$. The choice of the background metrics $h(t)$ is not important. In particular, we can choose the background metrics $h(t)$ to be independent of $t$. For each $t$, we denote by $\Delta_{g(t),h(t)}$ the Laplacian of a map from $(M,g(t))$ to $(M,h(t))$ (see \cite{Brendle2}, Definition 2.2). With this understood, we can define Ricci-DeTurck flow as follows:

\begin{definition}
Let $\tilde{g}(t)$ be a one-parameter family of metrics on $M$. We say that the metrics $\tilde{g}(t)$ evolve by the Ricci-DeTurck flow if 
\[\frac{\partial}{\partial t} \tilde{g}(t) = -2 \, \text{\rm Ric}_{\tilde{g}(t)} - \mathscr{L}_{\xi_t}(\tilde{g}(t)),\]
where the vector field $\xi_t$ is defined by $\xi_t := \Delta_{\tilde{g}(t),h(t)} \text{\rm id}$. 
\end{definition}

The evolution of the metric under the Ricci-DeTurck flow can be written in the form 
\[\frac{\partial}{\partial t} \tilde{g}_{ij} = \tilde{g}^{kl} \, \partial_k \partial_l \tilde{g}_{ij} + \text{\rm lower order terms.}\] 
Therefore, the Ricci-DeTurck flow is strictly parabolic, and admits a unique solution on a short time interval. 

In the next step, we show that the Ricci flow is equivalent to the Ricci-DeTurck flow in the sense that whenever we have a solution to one equation, we can convert it into a solution of the other.

To explain this, suppose first that we are given a solution $\tilde{g}(t)$ of the Ricci-DeTurck flow. Our goal is to produce a solution $g(t)$ of the Ricci flow. As above, we define $\xi_t := \Delta_{\tilde{g}(t),h(t)} \text{\rm id}$. We define a one-parameter family of diffeomorphisms $\varphi_t$ by $\frac{\partial}{\partial t} \varphi_t(p) = \xi_t |_{\varphi_t(p)}$ and $\varphi_0(p) = p$. Moreover, we define a one-parameter family of metrics $g(t)$ by $g(t) = \varphi_t^*(\tilde{g}(t))$. Then $g(t)$ is a solution of the Ricci flow. 

Conversely, suppose that we are given a solution $g(t)$ of the Ricci flow. Our goal is to produce a solution $\tilde{g}(t)$ of the Ricci-DeTurck flow. To that end, we solve the harmonic map heat flow $\frac{\partial}{\partial t} \varphi_t = \Delta_{g(t),h(t)} \varphi_t$ with initial condition $\varphi_0 = \text{\rm id}$. Moreover, we define a one-parameter family of metrics $\tilde{g}(t)$ by $\varphi_t^*(\tilde{g}(t)) = g(t)$. Then $\tilde{g}(t)$ is a solution of the Ricci-DeTurck flow. This shows that the Ricci flow is equivalent to the Ricci-DeTurck flow. 

A solution to the Ricci flow on a compact manifold can either be continued for all time, or else the curvature must blow up in finite time:

\begin{theorem}[R.~Hamilton \cite{Hamilton1}]
\label{curvature.unbounded.at.singularity}
Let $g_0$ be a Riemannian metric on a compact manifold $M$. Let $g(t)$, $t \in [0,T)$, denote the unique maximal solution to the Ricci flow with initial metric $g(0)=g_0$. If $T < \infty$, then the curvature of $g(t)$ is unbounded as $t \to T$.
\end{theorem}

A central problem is to understand the formation of singularities under the Ricci flow. To that end, it is often useful to consider a special class of solutions which move in a self-similar fashion. These are referred to as Ricci solitons:

\begin{definition}
Let $(M,g)$ be a Riemannian manifold, and let $f$ be a scalar function on $M$. We say that $(M,g,f)$ is a steady gradient Ricci soliton if $\text{\rm Ric} = D^2 f$. We say that $(M,g,f)$ is a shrinking gradient Ricci soliton if $\text{\rm Ric} = D^2 f + \mu g$ for some constant $\mu > 0$. We say that $(M,g,f)$ is an expanding gradient Ricci soliton if $\text{\rm Ric} = D^2 f + \mu g$ for some constant $\mu < 0$. 
\end{definition}

We next discuss the global behavior of the Ricci flow in dimension $2$. Hamilton \cite{Hamilton2} and Chow \cite{Chow1} showed that, for every initial metric on $S^2$, the Ricci flow shrinks to a point and becomes round after rescaling:

\begin{theorem}[R.~Hamilton \cite{Hamilton2}; B.~Chow \cite{Chow1}]
\label{2d.convergence}
Let $g_0$ be a Riemannian metric on $S^2$. Let $g(t)$, $t \in [0,T)$, denote the unique maximal solution to the Ricci flow with initial metric $g(0)=g_0$. Then $T < \infty$. Moreover, as $t \to T$, the rescaled metrics $\frac{1}{2(T-t)} \, g(t)$ converge in $C^\infty$ to a metric with constant Gaussian curvature $1$.
\end{theorem}

Theorem \ref{2d.convergence} was first proved by Hamilton \cite{Hamilton2} under the additional assumption that the initial metric $g_0$ has positive scalar curvature. This condition was later removed by Chow \cite{Chow1}. In the following, we sketch the main ideas in Hamilton's proof. Full details can be found in \cite{Hamilton2} or \cite{Brendle2}, Section 4. Given a metric $g$ on $S^2$ with positive scalar curvature, Hamilton defines the entropy $\mathcal{E}(g)$ by 
\begin{equation} 
\label{entropy}
\mathcal{E}(g) = \int_{S^2} R \, \log \Big ( \frac{AR}{8\pi} \Big ) \, d\mu, 
\end{equation} 
where $A$ denotes the area of $(S^2,g)$. The functional $\mathcal{E}(g)$ is invariant under scaling. By the Gauss-Bonnet theorem, $\int_{S^2} R \, d\mu = 8\pi$. Hence, it follows from Jensen's inequality that $\mathcal{E}(g)$ is nonnegative. Moreover, $\mathcal{E}(g)$ is strictly positive unless the scalar curvature of $(S^2,g)$ is constant. 

Hamilton's key insight is that the functional $\mathcal{E}(g)$ is monotone decreasing under the Ricci flow. From this, Hamilton deduced that the product $AR$ is uniformly bounded under the evolution. This implies that the flow converges to a shrinking gradient Ricci soliton, up to scaling. Finally, Hamilton showed that every shrinking gradient Ricci soliton on $S^2$ must have constant scalar curvature. This completes our discussion of Theorem \ref{2d.convergence}. \\

In the three-dimensional case, Hamilton \cite{Hamilton1} showed that an initial metric with positive Ricci curvature shrinks to a point in finite time and becomes round after rescaling.

\begin{theorem}[R.~Hamilton \cite{Hamilton1}]
\label{3d.convergence}
Let $g_0$ be a Riemannian metric on a three-manifold $M$ with positive Ricci curvature. Let $g(t)$, $t \in [0,T)$, denote the unique maximal solution to the Ricci flow with initial metric $g(0)=g_0$. Then $T < \infty$. Moreover, as $t \to T$, the rescaled metrics $\frac{1}{4(T-t)} \, g(t)$ converge in $C^\infty$ to a metric with constant sectional curvature $1$.
\end{theorem}

The proof of Theorem \ref{3d.convergence} is based on a pinching estimate for the eigenvalues of the Ricci tensor. To explain this, let $\lambda_1 \leq \lambda_2 \leq \lambda_3$ denote the eigenvalues of the tensor $R \, g_{ij} - 2 \, \text{\rm Ric}_{ij}$. With this understood, the scalar curvature is given by $\lambda_1+\lambda_2+\lambda_3$, and the eigenvalues of the Ricci tensor are given by $\frac{1}{2} \, (\lambda_2+\lambda_3),\frac{1}{2} \, (\lambda_3+\lambda_1),\frac{1}{2} \, (\lambda_1+\lambda_2)$. In particular, the positivity of the Ricci tensor is equivalent to the inequality $\lambda_1+\lambda_2 > 0$. Hamilton proved that 
\begin{equation} 
\frac{\partial}{\partial t} \lambda_1 \geq \Delta \lambda_1 + \lambda_1^2 + \lambda_2\lambda_3 
\end{equation} 
and 
\begin{equation} 
\frac{\partial}{\partial t} \lambda_3 \leq \Delta \lambda_3 + \lambda_3^2 + \lambda_1\lambda_2, 
\end{equation}
where both inequalities are understood in the barrier sense. In the special case when the initial metric has positive Ricci curvature, Hamilton proved a pinching estimate of the form $\lambda_3-\lambda_1 \leq C \, (\lambda_1+\lambda_2)^{1-\delta}$, where $\delta$ is a small positive constant depending on the initial data and $C$ is a large constant depending on the initial data. The proof of this estimate relies on the maximum principle.

Theorem \ref{3d.convergence} has opened up two major lines of research. On the one hand, it is of interest to prove similar convergence theorems in higher dimensions, under suitable assumptions on the curvature. This direction led to the proof of the Differentiable Sphere Theorem (see \cite{Brendle1},\cite{Brendle-Schoen}). On the other hand, it is important to understand the behavior of the Ricci flow in dimension $3$ for arbitrary initial metrics. In this case, the flow will develop more complicated types of singularities, including so-called neck-pinch singularities. In a series of breakthroughs, Perelman \cite{Perelman1},\cite{Perelman2} achieved a qualitative understanding of singularity formation in dimension $3$. This is sufficient for topological conclusions, such as the Poincar\'e conjecture.

In this survey, we will focus on issues related to singularity formation in dimension $3$. In Section \ref{ancient}, we will review the concept of an ancient solution, and explain its relevance for the analysis of singularities. We next discuss Perelman's noncollapsing theorem. Moreover, we describe examples of ancient solutions to the Ricci flow in low dimensions. In Section \ref{2d.classification}, we discuss the classification of ancient solutions in dimension $2$. In Section \ref{structure}, we review results due to Perelman \cite{Perelman1} concerning the structure of ancient $\kappa$-solutions in dimension $3$. These are ancient solutions which have bounded and nonnegative curvature and satisfy a noncollapsing condition. In Section \ref{3d.classification}, we discuss the classification of ancient $\kappa$-solutions in dimension $3$. In Section \ref{symmetry}, we describe a quantitative version of the fact that the Ricci flow preserves symmetry. In Section \ref{symmetry.improvement}, we discuss the Neck Improvement Theorem from \cite{Brendle4}. This theorem asserts that a neck becomes more symmetric under the evolution. Finally, in Sections \ref{soliton.asymptotics} and \ref{symmetry.of.solitons}, we give an alternative proof of the classification of noncollapsed steady gradient Ricci solitons in dimension $3$. This result was originally proved in \cite{Brendle3}; the proof given here relies on the Neck Improvement Theorem from \cite{Brendle4}.

\section{Ancient solutions and noncollapsing}

\label{ancient}

The notion of an ancient solution plays a fundamental role in understanding the formation of singularities in the Ricci flow. This concept was introduced by Hamilton \cite{Hamilton3}. 

\begin{definition}
An ancient solution to the Ricci flow is a solution which is defined on the time interval $(-\infty,T]$ for some $T$.
\end{definition}

The concept of an ancient solution to a parabolic PDE is analogous to the concept of an entire solution to an elliptic PDE.

Ancient solutions typically arise as blow-up limits at a singularity. In the Ricci flow, we are specifically interested in ancient solutions which satisfy a noncollapsing condition.

\begin{definition}[G.~Perelman \cite{Perelman1}]
\label{definition.of.noncollapsing}
An ancient solution to the Ricci flow in dimension $n$ is said to be $\kappa$-noncollapsed if $\text{\rm vol}_{g(t)}(B_{g(t)}(p,r)) \geq \kappa r^n$ whenever $\sup_{x \in B_{g(t)}(p,r)} R(x,t) \leq r^{-2}$.
\end{definition}

Definition \ref{definition.of.noncollapsing} is motivated by Perelman's noncollapsing theorem for the Ricci flow:

\begin{theorem}[G.~Perelman \cite{Perelman1}, Section 4]
\label{perelman.noncollapsing}
Let $M$ be a compact manifold of dimension $n$, and let $g(t)$, $t \in [0,T)$, be a solution to the Ricci flow, where $T<\infty$. Consider a sequence of times $t_j \to T$ and a bounded sequence of radii $r_j$. Finally, let $p_j$ be a sequence of points in $M$ such that 
\[r_j^2 \sup_{x \in B_{g(t_j)}(p_j,r_j)} R(x,t_j) < \infty.\] 
Then 
\[\liminf_{j \to \infty} r_j^{-n} \, \text{\rm vol}_{g(t_j)}(B_{g(t_j)}(p_j,r_j)) > 0.\]
\end{theorem} 

Theorem \ref{perelman.noncollapsing} is a consequence of Perelman's monotonicity formula for the $\mathcal{W}$-functional. In particular, Theorem \ref{perelman.noncollapsing} implies that every blow-up limit of the Ricci flow at a finite-time singularity must be $\kappa$-noncollapsed. 

In dimension $3$, the Hamilton-Ivey estimate gives a lower bound for the sectional curvature in terms of the scalar curvature:

\begin{theorem}[R.~Hamilton \cite{Hamilton3}; T.~Ivey \cite{Ivey}]
\label{hamilton.ivey}
Let $g(t)$, $t \in [0,T)$, be a solution to the Ricci flow on a compact three-manifold $M$. Let $\lambda_1$ denote the smallest eigenvalue of the tensor $R \, g_{ij} - 2 \, \text{\rm Ric}_{ij}$. Then $\lambda_1$ satisfies a pointwise inequality of the form $\lambda_1 \geq -f(R)$, where the function $f$ satisfies $\lim_{s \to \infty} \frac{f(s)}{s} = 0$.
\end{theorem}

Theorem \ref{hamilton.ivey} implies that every blow-up limit of the Ricci flow in dimension $3$ must have nonnegative sectional curvature. The proof of Theorem \ref{hamilton.ivey} relies on the maximum principle together with the evolution equation for the Ricci tensor. 

This motivates the following definition:

\begin{definition} 
An ancient $\kappa$-solution to the Ricci flow in dimension $n \in \{2,3\}$ is a complete, non-flat, $\kappa$-noncollapsed ancient solution with bounded and nonnegative curvature. 
\end{definition}

The notion of an ancient $\kappa$-solutions plays a key role in Perelman's theory. In particular, Perelman showed that, if a solution to the Ricci flow in dimension $3$ forms a singularity in finite time, then the high curvature regions can be approximated by ancient $\kappa$-solutions (see \cite{Perelman1}, Section 12).

In the remainder of this section, we describe some of the known examples of ancient solutions to the Ricci flow in dimension $2$ and $3$.

\begin{example}
Let $g_{S^2}$ denote the standard metric on $S^2$. Let us define a family of metrics $g(t)$ on $S^2$ by $g(t) = (-2t) \, g_{S^2}$ for $t \in (-\infty,0)$. This is an ancient solution to the Ricci flow which shrinks homothetically. It is $\kappa$-noncollapsed.
\end{example}

\begin{example}
Let us define a one-parameter family of conformal metrics on $\mathbb{R}^2$ by 
\[g_{ij}(t) = \frac{4}{e^t + |x|^2} \, \delta_{ij}\] 
for $t \in (-\infty,\infty)$. This gives a rotationally symmetric solution to the Ricci flow on $\mathbb{R}^2$, which moves by diffeomorphisms. It is referred to as the cigar soliton. The cigar soliton has positive curvature and opens up like a cylinder near infinity. The cigar soliton fails to be $\kappa$-noncollapsed.
\end{example}

\begin{example}
Let us define a one-parameter family of conformal metrics on $\mathbb{R}^2$ by 
\[g_{ij}(t) =  \frac{8 \, \sinh(-t)}{1 + 2 \cosh(-t) \, |x|^2 + |x|^4} \, \delta_{ij}\] 
for $t \in (-\infty,0)$. For each $t \in (-\infty,0)$, $g(t)$ extends to a smooth metric on $S^2$. This gives a rotationally symmetric solution to the Ricci flow on $S^2$. This is referred to as the King-Rosenau solution (cf. \cite{King},\cite{Rosenau}). The King-Rosenau solution is an ancient solution to the Ricci flow with positive curvature. The King-Rosenau solution fails to be $\kappa$-noncollapsed.
\end{example}

\begin{example} 
Let $g_{S^3}$ denote the standard metric on $S^3$. Let us define a family of metrics $g(t)$ on $S^2$ by $g(t) = (-4t) \, g_{S^3}$ for $t \in (-\infty,0)$. This is an ancient solution to the Ricci flow which shrinks homothetically. It is $\kappa$-noncollapsed.
\end{example}

\begin{example} 
Let again $g_{S^2}$ denote the standard metric on $S^2$. Let us define a family of metrics $g(t)$ on $S^2 \times \mathbb{R}$ by $g(t) = (-2t) \, g_{S^2} + dz \otimes dz$ for $t \in (-\infty,0)$. This is an ancient solution to the Ricci flow. It is $\kappa$-noncollapsed.
\end{example}

\begin{example}
Robert Bryant \cite{Bryant} has constructed a steady gradient Ricci soliton in dimension $3$ which is rotationally symmetric. This can be viewed as the three-dimensional analogue of the cigar soliton. The Bryant soliton has positive sectional curvature and opens up like a paraboloid near infinity. Unlike the cigar soliton, the Bryant soliton is $\kappa$-noncollapsed.
\end{example}

\begin{example}
Perelman has constructed an ancient solution to the Ricci flow on $S^3$ which is rotationally symmetric. This can be viewed as the three-dimensional analogue of the King-Rosenau solution. Perelman's ancient solution has positive sectional curvature. Unlike the King-Rosenau solution, Perelman's ancient solution is $\kappa$-noncollapsed. 

The asymptotics of Perelman's ancient solution are by now well understood; see \cite{Angenent-Brendle-Daskalopoulos-Sesum}.
\end{example}

\section{Classification of ancient solutions in dimension $2$}

\label{2d.classification}

In this section, we discuss the main classification results for ancient solution in dimension $2$. In \cite{Perelman1}, Perelman gave a classification of ancient $\kappa$-solutions in dimension $2$:

\begin{theorem}[G.~Perelman \cite{Perelman1}, Section 11]
\label{2D.noncollapsed.case}
Let $(M,g(t))$ be an ancient $\kappa$-solution in dimension $2$. Then $(M,g(t))$ is isometric to a family of shrinking spheres, or a $\mathbb{Z}_2$-quotient thereof.
\end{theorem}

Let us sketch Perelman's proof of Theorem \ref{2D.noncollapsed.case}. Suppose that $(M,g(t))$ is an ancient $\kappa$-solution in dimension $2$. After passing to a double cover if necessary, we may assume that $M$ is orientable. By Proposition 11.2 in \cite{Perelman1}, we can find a sequence of times $t_j \to -\infty$ and a sequence of points $p_j \in M$ with the following property: if we dilate the manifold $(M,g(t_j))$ around the point $p_j$ by the factor $(-t_j)^{-\frac{1}{2}}$, then the rescaled manifolds converge in the Cheeger-Gromov sense to a non-flat shrinking gradient Ricci soliton. Using Hamilton's classification of shrinking gradient Ricci solitons in dimension $2$ (see \cite{Hamilton2}), we conclude that the limiting manifold must be a round sphere. Since the limiting manifold is diffeomorphic to $S^2$, it follows that $M$ is diffeomorphic to $S^2$. 

We next consider Hamilton's entropy functional defined in (\ref{entropy}). Since the manifolds $(M,g(t_j))$ converge to a round sphere after rescaling, we know that $\mathcal{E}(g(t_j)) \to 0$ as $j \to \infty$. Moreover, it follows from Hamilton's work \cite{Hamilton2} that the function $t \mapsto \mathcal{E}(g(t))$ is monotone decreasing. Consequently, $\mathcal{E}(g(t)) \leq 0$ for each $t$. On the other hand, Jensen's inequality implies that $\mathcal{E}(g(t))$ is strictly positive unless $(M,g(t))$ has constant scalar curvature. Putting these facts together, we conclude that the scalar curvature of $(M,g(t))$ is constant for each $t$. This completes our sketch of the proof of Theorem \ref{2D.noncollapsed.case}. \\

Daskalopoulos, Hamilton, and \v Se\v sum were able to classify compact ancient solutions in dimension $2$ without noncollapsing assumptions: 

\begin{theorem}[P.~Daskalopoulos, R.~Hamilton, N.~\v Se\v sum \cite{Daskalopoulos-Hamilton-Sesum}]
\label{2D.collapsed.case}
Let $(M,g(t))$ be a compact, non-flat ancient solution to the Ricci flow in dimension $2$. Then, up to parabolic rescaling, translation in time, and diffeomorphisms, $(M,g(t))$ coincides with the family of shrinking spheres, the King-Rosenau solution, or a $\mathbb{Z}_2$-quotient of these.
\end{theorem}

One of the main ideas behind Theorem \ref{2D.collapsed.case} is to find a quantity to which the maximum principle can be applied and which vanishes on the King-Rosenau solution. To explain this, suppose that $(M,g(t))$ is a compact, non-flat ancient solution to the Ricci flow in dimension $2$. After passing to a double cover if necessary, we may assume that $M$ is orientable. Using the maximum principle, it is easy to see that $(M,g(t))$ has nonnegative scalar curvature for each $t$. Moreover, the strict maximum principle implies that the scalar curvature of $(M,g(t))$ is strictly positive. Since $M$ is compact and orientable, it follows that $M=S^2$. By the uniformization theorem, we may assume that the metrics $g(t)$ are conformal to the standard metric on $S^2$.

Using the stereographic projection, we can identify $\mathbb{R}^2$ with the complement of the north pole in $S^2$. Thus, we obtain a family of conformal metrics on $\mathbb{R}^2$ which evolve by the Ricci flow. We may write the evolving metric in the form $v^{-1} \, \delta_{ij}$, where $v$ satisfies the parabolic PDE 
\[\frac{\partial}{\partial t} v =  v \, \Delta v - |\nabla v|^2\] 
on $\mathbb{R}^2$, where $\nabla v$ and $\Delta v$ denote the gradient and Laplacian of $v$ with respect to the Euclidean metric on $\mathbb{R}^2$. Daskalopoulos, Hamilton, and \v Se\v sum consider the quantity 
\[Q = v \, \frac{\partial^3 v}{\partial z^3} \, \frac{\partial^3 v}{\partial \bar{z}^3},\] 
where $\frac{\partial}{\partial z} = \frac{1}{2} \, (\frac{\partial}{\partial x_1} - i \, \frac{\partial}{\partial x_2})$ and $\frac{\partial}{\partial \bar{z}} = \frac{1}{2} \, (\frac{\partial}{\partial x_1} + i \, \frac{\partial}{\partial x_2})$ denote the usual Wirtinger derivatives. A straightforward calculation shows that the quantity $Q$ is invariant under M\"obius transformations. Moreover, $Q$ satisfies the evolution equation 
\begin{align*} 
\frac{\partial}{\partial t} Q 
&= v \, \Delta Q - 16v^{-1} \, \Big ( v \, \frac{\partial^2 v}{\partial z \partial \bar{z}} - \frac{\partial v}{\partial z} \, \frac{\partial v}{\partial \bar{z}} \Big ) \, Q \\ 
&- 4 \, \Big ( v \, \frac{\partial^4 v}{\partial z^4} + 2 \, \frac{\partial v}{\partial z} \, \frac{\partial^3 v}{\partial z^3} \Big ) \, \Big ( v \, \frac{\partial^4 v}{\partial \bar{z}^4} + 2 \, \frac{\partial v}{\partial \bar{z}} \, \frac{\partial^3 v}{\partial \bar{z}^3} \Big ) \\ 
&- 4 \, \Big ( v \, \frac{\partial^4 v}{\partial z^3 \partial \bar{z}} - \frac{\partial v}{\partial \bar{z}} \, \frac{\partial^3 v}{\partial z^3} \Big ) \, \Big ( v \, \frac{\partial^4 v}{\partial \bar{z}^3 \partial z} - \frac{\partial v}{\partial z} \, \frac{\partial^3 v}{\partial \bar{z}^3} \Big ). 
\end{align*} 
The scalar curvature of the conformal metric $v^{-1} \, \delta_{ij}$ can be written in the form 
\[R = v^{-1} \, (v \, \Delta v - |\nabla v|^2) = 4v^{-1} \, \Big ( v \, \frac{\partial^2 v}{\partial z \partial \bar{z}} - \frac{\partial v}{\partial z} \, \frac{\partial v}{\partial \bar{z}} \Big ).\] 
This gives 
\[\frac{\partial}{\partial t} Q \leq v \, \Delta Q - 4RQ\] 
(compare \cite{Daskalopoulos-Hamilton-Sesum}, Section 5, or \cite{Chow2}). Here, $\Delta Q$ denotes the Laplacian of $Q$ with respect to the Euclidean metric on $\mathbb{R}^2$. The term $v \, \Delta Q$ can be interpreted as the Laplacian of $Q$ with respect to the evolving metric $v^{-1} \, \delta_{ij}$.

On the King-Rosenau solution, $v$ is a quadratic polynomial in $|x|^2$, with coefficients that depend on $t$. In particular, $\frac{\partial^3 v}{\partial z^3}$ vanishes identically on the King-Rosenau solution. Therefore, $Q$ vanishes identically on the King-Rosenau solution.

\section{Structure of ancient $\kappa$-solutions in dimension $3$}

\label{structure}

In \cite{Perelman1}, Perelman proved several fundamental results concerning the structure of ancient $\kappa$-solutions in dimension $3$. One of the central results is the following pointwise estimate for the covariant derivatives of the curvature tensor.

\begin{theorem}[G.~Perelman \cite{Perelman1}, Section 11]
\label{pointwise.derivative.estimate}
Let $(M,g(t))$, $t \in (-\infty,0]$, be an ancient $\kappa$-solution to the Ricci flow in dimension $3$. Let $m$ be a positive integer. Then the $m$-th order covariant derivatives of the curvature tensor satisfy the pointwise bound $|D^m \text{\rm Rm}| \leq C \, R^{\frac{m+2}{2}}$, where $C$ is a positive constant that depends only on $m$ and $\kappa$.
\end{theorem}

Another fundamental result in Perelman's work is the following longrange curvature estimate:

\begin{theorem}[G.~Perelman \cite{Perelman1}, Section 11]
\label{longrange.curvature.estimate}
Let $(M,g(t))$, $t \in (-\infty,0]$, be an ancient $\kappa$-solution to the Ricci flow in dimension $3$. Then there exists a function $\omega: [0,\infty) \to [0,\infty)$ (depending on $\kappa$) such that 
\[R(y,t) \leq R(x,t) \, \omega(R(x,t) \, d_{g(t)}(x,y)^2)\] 
for all $x,y \in M$ and all $t \leq 0$. 
\end{theorem}

Perelman's longrange curvature estimate is extremely useful, in that it allows Perelman to take limits of sequences of ancient $\kappa$-solutions. One important consequence is that the space of ancient $\kappa$-solutions is compact in the following sense: 

\begin{theorem}[G.~Perelman \cite{Perelman1}, Section 11]
\label{limits.of.ancient.solutions}
Let $(M^{(j)},g^{(j)}(t))$, $t \in (-\infty,0]$, be a sequence of ancient $\kappa$-solutions in dimension $3$. Moreover, suppose that $p_j \in M^{(j)}$ is a sequence of points such that $R(p_j,0) = 1$ for each $j$. Then, after passing to a subsequence if necessary, the flows $(M^{(j)},g^{(j)}(t),p_j)$ converge in the Cheeger-Gromov sense to a limit $(M^\infty,g^\infty(t))$, and this limit is again an ancient $\kappa$-solution.
\end{theorem}

\begin{corollary}[G.~Perelman \cite{Perelman1}, Section 11]
\label{limit.of.rescalings.of.an.ancient.solution}
Let $(M,g(t))$, $t \in (-\infty,0]$, be a noncompact ancient $\kappa$-solution to the Ricci flow in dimension $3$ with positive sectional curvature. Let us fix a point $p_0$ in $M$. Let $p_j$ be a sequence of points in $M$ such that $d_{g(0)}(p_0,p_j) \to \infty$, and let $r_j^{-2} := R(p_j,0)$. Let us dilate the flow around the point $(p_j,0)$ by the factor $r_j^{-1}$. Then, after passing to a subsequence if necessary, the rescaled flows converge in the Cheeger-Gromov sense to a family of shrinking cylinders.
\end{corollary}

Let us sketch how Corollary \ref{limit.of.rescalings.of.an.ancient.solution} follows from Theorem \ref{limits.of.ancient.solutions}. The longrange curvature estimate gives
\[R(p_0,0) \leq R(p_j,0) \, \omega(R(p_j,0) \, d_{g(0)}(p_0,p_j)^2)\] 
for each $j$. This implies 
\[R(p_j,0) \, d_{g(0)}(p_0,p_j)^2 \to \infty\] 
as $j \to \infty$. In other words, $r_j^{-1} \, d_{g(0)}(p_0,p_j) \to \infty$ as $j \to \infty$. We next consider the rescaled metrics $g^{(j)}(t) := r_j^{-2} \, g(r_j^2 t)$ for $t \in (-\infty,0]$. By Theorem \ref{limits.of.ancient.solutions}, the flows $(M,g^{(j)}(t),p_j)$ converge in the Cheeger-Gromov sense to an ancient $\kappa$-solution $(M^\infty,g^\infty(t))$. Since $r_j^{-1} \, d_{g(0)}(p_0,p_j) \to \infty$, the limiting flow $(M^\infty,g^\infty(t))$ must split off a line. Using Perelman's classification of ancient $\kappa$-solutions in dimension $2$ (see Theorem \ref{2D.noncollapsed.case}), it follows that the limiting flow $(M^\infty,g^\infty(t))$ must be a family of shrinking cylinders or a quotient thereof. On the other hand, $M$ is diffeomorphic to $\mathbb{R}^3$ by the soul theorem. In particular, $M$ does not contain an embedded $\mathbb{RP}^2$. This implies that $(M^\infty,g^\infty(t))$ cannot be a non-trivial quotient of the cylinder. This completes the sketch of the proof of Corollary \ref{limit.of.rescalings.of.an.ancient.solution}. \\

\begin{definition}
\label{evolving.neck}
Let $(M,g(t))$ be a solution to the Ricci flow in dimension $3$, and let $(\bar{x},\bar{t})$ be a point in space-time with $R(\bar{x},\bar{t}) = r^{-2}$. We say that $(\bar{x},\bar{t})$ lies at the center of an evolving $\varepsilon$-neck if, after rescaling by the factor $r^{-1}$, the parabolic neighborhood $B_{g(\bar{t})}(\bar{x},\varepsilon^{-1} r) \times [\bar{t}-\varepsilon^{-1} r^2,\bar{t}]$ is $\varepsilon$-close in $C^{[\varepsilon^{-1}]}$ to a family of shrinking cylinders.
\end{definition}

The notion of a neck was introduced in Hamilton's work \cite{Hamilton4}. In particular, Hamilton showed that a neck admits a canonical foliation by constant mean curvature (CMC) spheres.

Corollary \ref{limit.of.rescalings.of.an.ancient.solution} implies the following structure theorem for noncompact ancient $\kappa$-solutions:

\begin{corollary}[G.~Perelman \cite{Perelman1}, Section 11]
\label{structure.theorem}
Let $(M,g(t))$, $t \in (-\infty,0]$, be a noncompact ancient $\kappa$-solution to the Ricci flow in dimension $3$ with positive sectional curvature. Moreover, let $\varepsilon$ be a positive real number, and let $M_\varepsilon$ denote the set of all points $x \in M$ with the property that $(x,0)$ does not lie at the center of an evolving $\varepsilon$-neck. Then $M_\varepsilon$ has finite diameter. Moreover, $\sup_{x \in M_\varepsilon} R(x,0) \leq C(\kappa,\varepsilon) \, \inf_{x \in M_\varepsilon} R(x,0)$ and $\sup_{x \in M_\varepsilon} R(x,0) \leq C(\kappa,\varepsilon) \, \text{\rm diam}_{g(0)}(M_\varepsilon)^{-2}$.
\end{corollary}

\section{Classification of ancient $\kappa$-solutions in dimension $3$} 

\label{3d.classification}

We now turn to the classification of ancient $\kappa$-solutions in dimension $3$. The first major step was the classification of noncollapsed steady gradient Ricci solitons in dimension $3$.

\begin{theorem}[S.~Brendle \cite{Brendle3}]
\label{soliton.classification}
Let $(M,g)$ be a three-dimensional complete steady gradient Ricci soliton which is non-flat and $\kappa$-noncollapsed. Then $(M,g)$ is rotationally symmetric, and is therefore isometric to the Bryant soliton up to scaling. 
\end{theorem}

More recently, we classified all noncompact ancient $\kappa$-solutions in dimension $3$: 

\begin{theorem}[S.~Brendle \cite{Brendle4}]
\label{noncompact.case}
Assume that $(M,g(t))$ is a noncompact ancient $\kappa$-solution of dimension $3$. Then either $(M,g(t))$ is isometric to a family of shrinking cylinders (or a quotient thereof), or $(M,g(t))$ is isometric to the Bryant soliton up to scaling. 
\end{theorem}

Theorem \ref{noncompact.case} confirms a conjecture of Perelman \cite{Perelman1}.

The proof of Theorem \ref{noncompact.case} consists of two main steps. In the first step, we classify noncompact ancient $\kappa$-solutions with rotational symmetry. To do that, we need precise asymptotic estimates for such solutions. In the second step, we show that every noncompact ancient $\kappa$-solution is rotationally symmetric. This second step uses the classification of steady gradient Ricci solitons in Theorem \ref{soliton.classification}, as well as the classification of ancient $\kappa$-solutions with rotational symmetry. Another crucial ingredient is the Neck Improvement Theorem which asserts that a neck tends to get more symmetric as it evolves under the Ricci flow. We will discuss the Neck Improvement Theorem in Section \ref{symmetry.improvement} below.

The following theorem is the counterpart of Theorem \ref{noncompact.case} in the compact case:

\begin{theorem}[S.~Brendle, P.~Daskalopoulos, N.~\v Se\v sum \cite{Brendle-Daskalopoulos-Sesum}]
\label{compact.case}
Assume that $(M,g(t))$ is a compact ancient $\kappa$-solution of dimension $3$. Then, up to parabolic rescaling, translation in time, and diffeomorphisms, $(M,g(t))$ is either a family of shrinking spheres or Perelman's ancient solution or a quotient of these.
\end{theorem}

The proof of Theorem \ref{compact.case} again requires two main steps. In the first step, we show that every compact ancient $\kappa$-solution is rotationally symmetric. This step uses the classification of noncompact ancient $\kappa$-solutions in Theorem \ref{noncompact.case}, together with the Neck Improvement Theorem. In a second step, we classify compact ancient $\kappa$-solutions with rotational symmetry. To do that, we need to understand the asymptotic behavior of such solutions. These asymptotic estimates are established in \cite{Angenent-Brendle-Daskalopoulos-Sesum}.

Similar classification results exist for convex, noncollapsed ancient solutions to mean curvature flow in $\mathbb{R}^3$. We refer to \cite{Brendle-Choi} for the classification in the noncompact case, and to \cite{Angenent-Daskalopoulos-Sesum} for the classification in the compact case.

\section{Preservation of symmetry under the Ricci flow}

\label{symmetry}

Let $g(t)$, $t \in [0,T)$, be a solution to the Ricci flow on a compact manifold $M$. It follows from Hamilton's short time uniqueness theorem that the Ricci flow preserves symmetry. More precisely, every isometry of $(M,g(0))$ is also an isometry of $(M,g(t))$ for each $t \geq 0$. Consequently, every Killing vector field of $(M,g(0))$ is also a Killing vector field of $(M,g(t))$ for each $t \geq 0$.

In this section, we describe a quantitative version of this principle. We begin with a definition: 

\begin{definition}
Let $h$ be a symmetric $(0,2)$-tensor. The Lichnerowicz Laplacian of $h$ is defined as 
\begin{equation} 
\Delta_L h_{ik} := \Delta h_{ik} + 2 R_{ijkl} h^{jl} - \text{\rm Ric}_i^l h_{kl} - \text{\rm Ric}_k^l h_{il}. 
\end{equation} 
\end{definition}

The Lichnerowicz Laplacian arises naturally in the study of Einstein metrics (see \cite{Besse}, equation (1.180b)). It also comes up in connection with the evolution equation for the Ricci tensor under the Ricci flow. Indeed, if $(M,g(t))$ is a solution to the Ricci flow, then the Ricci tensor satisfies the evolution equation 
\[\frac{\partial}{\partial t} \text{\rm Ric}_{g(t)} = \Delta_{L,g(t)} \text{\rm Ric}_{g(t)}\] 
(see \cite{Hamilton1}, Corollary 7.3). 

The following results play a key role in our analysis:

\begin{proposition}
\label{divergence.of.Lie.derivative}
Let $(M,g)$ be a Riemannian manifold. Let $V$ be a smooth vector field on $M$, and let $h := \mathscr{L}_V(g)$. Then 
\[\Delta V + \text{\rm Ric}(V) = \text{\rm div} \, h - \frac{1}{2} \, \nabla(\text{\rm tr} \, h).\] 
\end{proposition}

\begin{proposition}[S.~Brendle \cite{Brendle4}, Section 5]
\label{evolution.of.Lie.derivative}
Let $(M,g(t))$ be a solution of the Ricci flow. Let $V(t)$ be a time-dependent vector field such that 
\begin{equation}
\label{pde.for.vector.field} 
\frac{\partial}{\partial t} V(t) = \Delta_{g(t)} V(t) + \text{\rm Ric}_{g(t)}(V(t)), 
\end{equation} 
and let $h(t) := \mathscr{L}_{V(t)}(g(t))$. Then the tensor $h(t)$ satisfies the evolution equation
\begin{equation}
\label{pde.for.h} 
\frac{\partial}{\partial t} h(t) = \Delta_{L,g(t)} h(t). 
\end{equation}
\end{proposition}

Proposition \ref{evolution.of.Lie.derivative} has a natural geometric interpretation in terms of the linearized Ricci-DeTurck flow. To explain this, let us fix a solution $(M,g(t))$ of the Ricci flow. Suppose that $\varphi_t$ is a one-parameter family of diffeomorphisms which solve the harmonic map heat flow with respect to the background metrics $g(t)$; that is, 
\begin{equation} 
\label{harmonic.map.heat.flow}
\frac{\partial}{\partial t} \varphi_t = \Delta_{g(t),g(t)} \varphi_t. 
\end{equation}
Let us define a flow of metrics $\tilde{g}(t)$ by $\varphi_t^*(\tilde{g}(t)) = g(t)$. Then the metrics $\tilde{g}(t)$ solve the Ricci-DeTurck flow with respect to the background metrics $g(t)$. More precisely,
\begin{equation} 
\label{ricci.deturck.flow} 
\frac{\partial}{\partial t} \tilde{g}(t) = -2 \, \text{\rm Ric}_{\tilde{g}(t)} - \mathscr{L}_{\xi_t}(\tilde{g}(t)), 
\end{equation}
where $\xi_t = \Delta_{\tilde{g}(t),g(t)} \text{\rm id}$. Clearly, $\varphi_t := \text{\rm id}$ is a solution of (\ref{harmonic.map.heat.flow}), and $\tilde{g}(t) := g(t)$ is a solution of (\ref{ricci.deturck.flow}). 

We now linearize the equations (\ref{harmonic.map.heat.flow}) and (\ref{ricci.deturck.flow}) around $\varphi_t = \text{\rm id}$ and $\tilde{g}(t) = g(t)$, respectively. Linearizing the harmonic map heat flow (\ref{harmonic.map.heat.flow}) around the identity, we obtain the equation 
\[\frac{\partial}{\partial t} V(t) = \Delta_{g(t)} V(t) + \text{\rm Ric}_{g(t)}(V(t))\] 
for a vector field $V$ (see \cite{Eells-Lemaire}, p.~11). Linearizing the Ricci-DeTurck flow (\ref{ricci.deturck.flow}) around $g(t)$ leads to the parabolic Lichnerowicz equation 
\[\frac{\partial}{\partial t} h(t) = \Delta_{L,g(t)} h(t).\] 
This completes our discussion of Proposition \ref{evolution.of.Lie.derivative}. \\

On a steady gradient Ricci soliton, Proposition \ref{evolution.of.Lie.derivative} takes the following form:

\begin{corollary}[S.~Brendle \cite{Brendle3}]
\label{evolution.of.Lie.derivative.soliton.version}
Let $(M,g,f)$ be a steady gradient Ricci soliton, and let $X := \nabla f$. Let $V$ be a vector field satisfing 
\begin{equation} 
\label{pde.for.vector.field.soliton.version} 
\Delta V + D_X V = 0, 
\end{equation} 
and let $h := \mathscr{L}_V(g)$. Then the tensor $h$ satisfies the equation 
\begin{equation} 
\label{pde.for.h.soliton.version} 
\Delta_L h + \mathscr{L}_X(h) = 0. 
\end{equation}
\end{corollary}

Let us sketch how Corollary \ref{evolution.of.Lie.derivative.soliton.version} follows from Theorem \ref{evolution.of.Lie.derivative}. On a steady gradient Ricci soliton, the time derivative $\frac{\partial}{\partial t}$ reduces to a Lie derivative $-\mathscr{L}_X$. More precisely, let $(M,g,f)$ be a steady gradient Ricci soliton, let $X := \nabla f$, and let $\Phi_t$ denote the flow generated by the vector field $-X$. Suppose that $V$ satisfies $\Delta V + D_X V = 0$, and let $h := \mathscr{L}_V(g)$. Using the identity $D_V X = \text{\rm Ric}(V)$, we obtain $\Delta V + \mathscr{L}_X V + \text{\rm Ric}(V) = 0$. Consequently, the vector fields $\Phi_t^*(V)$ satisfy the parabolic PDE (\ref{pde.for.vector.field}) on the evolving background $(M,\Phi_t^*(g))$. By Theorem \ref{evolution.of.Lie.derivative}, the tensors $\Phi_t^*(h)$ satisfy the parabolic PDE (\ref{pde.for.h}) on the evolving background $(M,\Phi_t^*(g))$. This implies $\Delta_L h + \mathscr{L}_X(h) = 0$.

In order to apply Proposition \ref{evolution.of.Lie.derivative} in practice, we need estimates for solutions of the parabolic Lichnerowicz equation. In dimension $3$, this can be accomplished by applying the maximum principle to the quantity $\frac{|h|^2}{R^2}$:

\begin{proposition}[G.~Anderson, B.~Chow \cite{Anderson-Chow}]
\label{anderson.chow.estimate}
Let $(M,g(t))$ be a solution to the Ricci flow in dimension $3$ with positive scalar curvature. Let $h$ be a solution of the parabolic Lichnerowicz equation $\frac{\partial}{\partial t} h(t) = \Delta_{L,g(t)} h(t)$. Then 
\[\frac{\partial}{\partial t} \Big ( \frac{|h|^2}{R^2} \Big ) \leq \Delta \Big ( \frac{|h|^2}{R^2} \Big ) + \frac{2}{R} \, \Big \langle \nabla R,\nabla \Big ( \frac{|h|^2}{R^2} \Big ) \Big \rangle.\] 
\end{proposition}

On a steady gradient Ricci soliton, Proposition \ref{anderson.chow.estimate} takes the following form:

\begin{corollary}[G.~Anderson, B.~Chow \cite{Anderson-Chow}]
\label{anderson.chow.estimate.soliton.version}
Let $(M,g,f)$ be a steady gradient Ricci soliton in dimension $3$ with positive scalar curvature, and let $X := \nabla f$. Let $h$ be a solution of the equation $\Delta_L h + \mathscr{L}_X(h) = 0$. Then 
\[\Delta \Big ( \frac{|h|^2}{R^2} \Big ) + \Big \langle X,\nabla \Big ( \frac{|h|^2}{R^2} \Big ) \Big \rangle + \frac{2}{R} \, \Big \langle \nabla R,\nabla \Big ( \frac{|h|^2}{R^2} \Big ) \Big \rangle \geq 0.\] 
\end{corollary}

\section{Improvement of symmetry on a neck}

\label{symmetry.improvement}

In view of Proposition \ref{evolution.of.Lie.derivative}, it is important to understand the parabolic Lichnerowicz equation (\ref{pde.for.h}) on a Ricci flow background. As a starting point, we consider the special case when the background is given by a family of shrinking cylinders. To fix notation, we define a family of metrics $\bar{g}(t)$, $t \in (-\infty,0)$ on $S^2 \times \mathbb{R}$ by  
\[\bar{g}(t) = (-2t) \, g_{S^2} + dz \otimes dz, \quad t \in (-\infty,0).\] 
Clearly, the metrics $\bar{g}(t)$, $t \in (-\infty,0)$, evolve by the Ricci flow. 

\begin{proposition}[S.~Brendle \cite{Brendle4}, Section 6]
\label{parabolic.Lichnerowicz.equation.on.cylinders}
Let $(S^2 \times \mathbb{R},\bar{g}(t))$ denote the family of shrinking cylinders. Let $L$ be a large real number. Let $h(t)$ be a solution of the parabolic Lichnerowicz equation $\frac{\partial}{\partial t} h(t) = \Delta_{L,\bar{g}(t)} h(t)$ which is defined on $S^2 \times [-\frac{L}{2},\frac{L}{2}]$ and for $t \in [-\frac{L}{2},-1]$. Assume that $|h(t)|_{\bar{g}(t)} \leq 1$ for $t \in [-\frac{L}{2},-\frac{L}{4}]$, and $|h(t)|_{\bar{g}(t)} \leq L^{10}$ for $t \in [-\frac{L}{4},-1]$. Then we can find a rotationally invariant tensor of the form $\bar{\omega}(z,t) \, g_{S^2} + \bar{\beta}(z,t) \, dz \otimes dz$ and a scalar function $\psi: S^2 \to \mathbb{R}$ such that $\psi$ lies in the span of the first spherical harmonics on $S^2$ and 
\[|h(t) - \bar{\omega}(z,t) \, g_{S^2} - \bar{\beta}(z,t) \, dz \otimes dz - (-t) \, \psi \, g_{S^2}|_{\bar{g}(t)} \leq C \, L^{-\frac{1}{2}}\] 
on $S^2 \times [-1000,1000]$ and for $t \in [-1000,-1]$. Here, $C$ is a constant which does not depend on $L$.
\end{proposition}

Proposition \ref{parabolic.Lichnerowicz.equation.on.cylinders} asserts that, given sufficient time to evolve, a solution of the parabolic Lichnerowicz equation can be approximated by a sum of a rotationally invariant tensor and a tensor of the form $(-t) \, \psi \, g_{S^2}$, where $\psi: S^2 \to \mathbb{R}$ lies in the span of the first spherical harmonics on $S^2$. The tensor $(-t) \, \psi \, g_{S^2}$ can be written as a Lie derivative of the metric along a vector field. To see this, let us define a vector field $\xi$ on $S^2$ by $g_{S^2}(\xi,\cdot) = -\frac{1}{4} \, d\psi$. Since $\psi$ lies in the span of the first spherical harmonics on $S^2$, we obtain $\mathscr{L}_\xi(g_{S^2}) = \frac{1}{2} \, \psi \, g_{S^2}$, and consequently $\mathscr{L}_\xi(\bar{g}(t)) = (-t) \, \psi \, g_{S^2}$.

To prove Proposition \ref{parabolic.Lichnerowicz.equation.on.cylinders}, we decompose the tensor $h(t)$ into components, and perform a mode decomposition in spherical harmonics. This leads to a system of linear heat equations in one space dimension. 

Using Proposition \ref{parabolic.Lichnerowicz.equation.on.cylinders}, we can show that a neck becomes more symmetric as it evolves under the Ricci flow. To state this result, we need a quantitative notion of $\varepsilon$-symmetry:

\begin{definition}[S.~Brendle \cite{Brendle4}, Section 8]
\label{symmetry.of.necks}
Let $(M,g(t))$ be a solution to the Ricci flow in dimension $3$, and let $(\bar{x},\bar{t})$ be a point in space-time with $R(\bar{x},\bar{t}) = r^{-2}$. We assume that $(\bar{x},\bar{t})$ lies at the center of an evolving $\varepsilon _0$-neck for some small positive number $\varepsilon_0$. We say that $(\bar{x},\bar{t})$ is $\varepsilon$-symmetric if there exist smooth, time-independent vector fields $U^{(1)},U^{(2)},U^{(3)}$ which are defined on an open set containing $\bar{B}_{g(\bar{t})}(\bar{x},100r)$ and satisfy the following conditions: 
\begin{itemize}
\item $\sup_{\bar{B}_{g(\bar{t})}(\bar{x},100r) \times [\bar{t}-100r^2,\bar{t}]} \sum_{l=0}^2 \sum_{a=1}^3 r^{2l} \, |D^l(\mathscr{L}_{U^{(a)}}(g(t)))|^2 \leq \varepsilon^2$.
\item If $t \in [\bar{t}-100r^2,\bar{t}]$ and $\Sigma \subset \bar{B}_{g(\bar{t})}(\bar{x},100r)$ is a leaf of the CMC foliation of $(M,g(t))$, then $\sup_\Sigma \sum_{a=1}^3 r^{-2} \, |\langle U^{(a)},\nu \rangle|^2 \leq \varepsilon^2$, where $\nu$ denotes the unit normal vector to $\Sigma$ in $(M,g(t))$.
\item If $t \in [\bar{t}-100r^2,\bar{t}]$ and $\Sigma \subset \bar{B}_{g(\bar{t})}(\bar{x},100r)$ is a leaf of the CMC foliation of $(M,g(t))$, then 
\[\sum_{a,b=1}^3 \bigg | \delta_{ab} - \text{\rm area}_{g(t)}(\Sigma)^{-2} \int_\Sigma \langle U^{(a)},U^{(b)} \rangle_{g(t)} \, d\mu_{g(t)} \bigg |^2 \leq \varepsilon^2.\]
\end{itemize}
\end{definition}

With this understood, we can now state the Neck Improvement Theorem from \cite{Brendle4}:

\begin{theorem}[S.~Brendle \cite{Brendle4}, Section 8]
\label{neck.improvement.theorem}
We can find a large constant $L$ and small positive constant $\varepsilon_1$ such that the following holds. Let $(M,g(t))$ be a solution of the Ricci flow in dimension $3$, and let $(x_0,t_0)$ be a point in space-time which lies at the center of an evolving $\varepsilon_1$-neck and satisfies $R(x_0,t_0) = r^{-2}$. Moreover, we assume that every point in the parabolic neighborhood $B_{g(t_0)}(x_0,Lr) \times [t_0-Lr^2,t_0)$ is $\varepsilon$-symmetric, where $\varepsilon \leq \varepsilon_1$. Then the point $(x_0,t_0)$ is $\frac{\varepsilon}{2}$-symmetric.
\end{theorem}

\section{Asymptotic behavior of noncollapsed steady gradient Ricci solitons in dimension $3$}

\label{soliton.asymptotics}

Let $(M,g,f)$ be a non-flat steady gradient Ricci soliton in dimension $n$, so that $\text{\rm Ric} = D^2 f$. For abbreviation, let $X := \nabla f$. Throughout this section, we fix an arbitrary point $p \in M$. 

\begin{lemma}
\label{barrier}
Given any point $\bar{x} \in M$, we can find a smooth function $\rho$ which is defined in an open neighborhood of $\bar{x}$ and satisfies the following conditions: 
\begin{itemize}
\item $\rho(x) \geq d(p,x)^2$ in an open neighborhood of the point $\bar{x}$.
\item $\rho(\bar{x}) = d(p,\bar{x})^2$.
\item $|\nabla \rho|^2 = 4\rho$ at the point $\bar{x}$.
\item $\Delta \rho + \langle X,\nabla \rho \rangle \leq N_0 + N_1 \sqrt{\rho}$ at the point $\bar{x}$.
\end{itemize}
Here, $N_0$ and $N_1$ are uniform constants which do not depend on $\bar{x}$.
\end{lemma}

\textbf{Proof.} 
Let us fix a positive real number $r_0$ such that $r_0$ is strictly smaller than the injectivity radius at $p$. 

We first consider the case $\bar{x} \in B(p,r_0)$. In this case, we define a smooth function $\rho: B(p,r_0) \to \mathbb{R}$ by $\rho(x) := d(p,x)^2$. It is easy to see that $|\nabla \rho|^2 = 4\rho$ at each point in $B(p,r_0)$. Moreover, at each point in $B(p,r_0)$, we have $\Delta \rho + \langle X,\nabla \rho \rangle \leq C$ for some uniform constant $C$.

In the next step, we consider the case $\bar{x} \in M \setminus B(p,r_0)$. Let $l := d(p,\bar{x}) \geq r_0$. Moreover, let $\gamma: [0,l] \to M$ be a unit-speed geodesic with $\gamma(0)=p$ and $\gamma(l)=\bar{x}$. Finally, let $\chi: [0,\infty) \to [0,\infty)$ be a smooth cutoff function such that $\chi=0$ on the interval $[0,\frac{r_0}{2}]$ and $\chi=1$ on the interval $[r_0,\infty)$. 

Let us fix a positive real number $\bar{r}$ such that $\bar{r}$ is strictly smaller than the injectivity radius at $\bar{x}$. We define a smooth function $\rho: B(\bar{x},\bar{r}) \to \mathbb{R}$ as follows. Given a point $x \in B(\bar{x},\bar{r})$, there exists a unique vector $w \in T_{\bar{x}} M$ such that $|w| < \bar{r}$ and $x = \exp_{\bar{x}}(w)$. We denote by $W$ the unique parallel vector field along $\gamma$ satisfying $W(l)=w$. We then define $\sqrt{\rho(x)}$ to be the length of the curve 
\[s \mapsto \exp_{\gamma(s)}(\chi(s) \, W(s)), \quad s \in [0,l].\] 
Clearly, 
\[\sqrt{\rho(x)} \geq d \big ( \exp_{\gamma(0)}(\chi(0) \, W(0)),\exp_{\gamma(l)}(\chi(l) \, W(l)) \big ) = d(p,\exp_{\bar{x}}(w)) = d(p,x).\] 
Moreover, in the special case when $x=\bar{x}$ and $w=0$, we obtain 
\[\sqrt{\rho(\bar{x})} = l = d(p,\bar{x}).\]
The formula for the first variation of arclength implies that $\nabla \sqrt{\rho} = \gamma'(l)$ at the point $\bar{x}$. In particular, $|\nabla \sqrt{\rho}|^2 = 1$ at the point $\bar{x}$. Consequently, $|\nabla \rho|^2 = 4\rho$ at the point $\bar{x}$.

Using the formula for the second variation of arclength, we obtain 
\begin{align*} 
(D^2 \sqrt{\rho})_{\bar{x}}(w,w) 
&= \int_0^l \chi'(s)^2 \, \langle W(s),W(s) \rangle \, ds - \int_0^l \chi'(s)^2 \, \langle \gamma'(s),W(s) \rangle^2 \, ds \\ 
&- \int_0^l \chi(s)^2 \, R(\gamma'(s),W(s),\gamma'(s),W(s)) \, ds 
\end{align*} 
for every vector $w \in T_{\bar{x}} M$, where $W$ denotes the unique parallel vector field along $\gamma$ satisfying $W(l)=w$. Taking the trace over $w$ gives 
\[\Delta \sqrt{\rho}(\bar{x}) = (n-1) \int_0^l \chi'(s)^2 \, ds - \int_0^l \chi(s)^2 \, \text{\rm Ric}(\gamma'(s),\gamma'(s)) \, ds.\] 
Since $\chi=1$ on the interval $[r_0,l]$, we obtain 
\[\Delta \sqrt{\rho}(\bar{x}) \leq -\int_0^l \text{\rm Ric}(\gamma'(s),\gamma'(s)) \, ds + C,\] 
where $C$ denotes a uniform constant that does not depend on $\bar{x}$. 

On the other hand, using the identity $D^2 f = \text{\rm Ric}$, we obtain 
\[\frac{d}{ds} \langle \nabla f(\gamma(s)),\gamma'(s) \rangle = \text{\rm Ric}(\gamma'(s),\gamma'(s)).\] 
Integrating this identity over $s \in [0,l]$ gives 
\[\langle \nabla f(\gamma(l)),\gamma'(l) \rangle \leq \int_0^l \text{\rm Ric}(\gamma'(s),\gamma'(s)) \, ds + C,\] 
where $C$ is a uniform constant that does not depend on $\bar{x}$. Since $\gamma(l) = \bar{x}$ and $\gamma'(l) = \nabla \sqrt{\rho}(\bar{x})$, we obtain 
\[\langle \nabla f(\bar{x}),\nabla \sqrt{\rho}(\bar{x}) \rangle \leq \int_0^l \text{\rm Ric}(\gamma'(s),\gamma'(s)) \, ds + C,\] 
where $C$ is a uniform constant that does not depend on $\bar{x}$. Putting these facts together, we conclude that 
\[\Delta \sqrt{\rho}(\bar{x}) + \langle \nabla f(\bar{x}),\nabla \sqrt{\rho}(\bar{x}) \rangle \leq C,\] 
where $C$ is a uniform constant that does not depend on $\bar{x}$. This finally implies 
\[\Delta \rho(\bar{x}) + \langle \nabla f(\bar{x}),\nabla \rho(\bar{x}) \rangle \leq 2 + C \sqrt{\rho},\] 
where $C$ is a uniform constant that does not depend on $\bar{x}$. This completes the proof of Lemma \ref{barrier}. \\

\begin{proposition}[B.L.~Chen \cite{Chen}]
\label{nonnegative.scalar.curvature}
The manifold $(M,g)$ has nonnegative scalar curvature.
\end{proposition}

\textbf{Proof.} 
As above, we fix an arbitrary point $p \in M$. Let $N_0$ and $N_1$ denote the constants in Lemma \ref{barrier}. Let us fix a radius $r > 0$. We define a continuous function $u: B(p,r) \to \mathbb{R}$ by 
\[u(x) := R(x) + 2n(12+N_0+N_1 r) r^2 (r^2-d(p,x)^2)^{-2}\] 
for $x \in B(p,r)$. We claim that $u(x) \geq 0$ for all $x \in B(p,r)$. To prove this, we argue by contradiction. Let $\bar{x}$ be a point in $B(p,r)$ where the function $u$ attains its minimum, and suppose that $u(\bar{x}) < 0$. The evolution equation for the scalar curvature implies 
\[\Delta R + \langle X,\nabla R \rangle = -2 \, |\text{\rm Ric}|^2.\] 
By Lemma \ref{barrier}, we can find a smooth function $\rho$ which is defined in an open neighborhood of $\bar{x}$ and satisfies the following conditions: 
\begin{itemize}
\item $\rho(x) \geq d(p,x)^2$ in an open neighborhood of the point $\bar{x}$.
\item $\rho(\bar{x}) = d(p,\bar{x})^2$.
\item $|\nabla \rho|^2 = 4\rho$ at the point $\bar{x}$.
\item $\Delta \rho + \langle X,\nabla \rho \rangle \leq N_0 + N_1 \sqrt{\rho}$ at the point $\bar{x}$.
\end{itemize}
Then 
\[R(x) + 2n(12+N_0+N_1 r) r^2 (r^2-\rho(x))^{-2} \geq u(x) \geq u(\bar{x})\] 
in an open neighborhood of $\bar{x}$, with equality at the point $\bar{x}$. Consequently, the function $R + 2n(12+N_0+ N_1r) r^2 (r^2-\rho)^{-2}$ attains a local minimum at the point $\bar{x}$. Thus, we conclude that 
\begin{align*} 
0 &\leq \Delta R + \langle X,\nabla R \rangle \\ 
&+ 2n(12+N_0+ N_1r) r^2 \, \big [ \Delta((r^2-\rho)^{-2}) + \langle X,\nabla((r^2-\rho)^{-2}) \rangle \big ] \\ 
&= \Delta R + \langle X,\nabla R \rangle \\ 
&+ 4n(12+N_0+N_1 r) r^2 (r^2-\rho)^{-3} \, (\Delta \rho + \langle X,\nabla \rho \rangle) \\ 
&+ 12n(12+N_0+N_1 r) r^2 (r^2-\rho)^{-4} \, |\nabla \rho|^2 \\ 
&\leq -2 \, |\text{\rm Ric}|^2 + 4n(12+N_0+N_1 r) (N_0+N_1 \sqrt{\rho}) r^2 (r^2-\rho)^{-3} \\ 
&+ 48n(12+N_0+N_1 r) \rho r^2 (r^2-\rho)^{-4} 
\end{align*} 
at the point $\bar{x}$. Since $u \leq 0$ at the point $\bar{x}$, we know that 
\[-R \geq 2n(12+N_0+N_1 r) r^2 (r^2-\rho)^{-2}\] 
at the point $\bar{x}$. This implies 
\[n \, |\text{\rm Ric}|^2 \geq R^2 \geq 4n^2(12+N_0+N_1 r)^2 r^4 (r^2-\rho)^{-4}\] 
at the point $\bar{x}$. Putting these facts together, we obtain 
\begin{align*} 
0 &\leq -2 \, |\text{\rm Ric}|^2 + 4n(12+N_0+N_1 r) (N_0+N_1 \sqrt{\rho}) r^2 (r^2-\rho)^{-3} \\ 
&+ 48n(12+N_0+N_1 r) \rho r^2 (r^2-\rho)^{-4} \\ 
&\leq - 8n(12+N_0+N_1 r)^2 r^4 (r^2-\rho)^{-4} \\ 
&+ 4n(12+N_0+N_1 r) (N_0+N_1 r) r^4 (r^2-\rho)^{-4} \\ 
&+ 48n(12+N_0+N_1 r) r^4 (r^2-\rho)^{-4} \\ 
&= -4n(12+N_0+N_1 r)^2 r^4 (r^2-\rho)^{-4}
\end{align*} 
at the point $\bar{x}$. This is a contradiction. 

Thus, we conclude that 
\[R(x) + 2n(12+N_0+N_1 r) r^2 (r^2-d(p,x)^2)^{-2} \geq 0\] 
for all $x \in B(p,r)$. Sending $r \to \infty$ gives $R(x) \geq 0$ for each point $x \in M$. This completes the proof of Proposition \ref{nonnegative.scalar.curvature}. \\

\begin{corollary}
\label{bound.for.gradient.of.f}
There exists a large constant $C$ such that $|\nabla f| \leq C$ at each point on $M$.
\end{corollary}

\textbf{Proof.} 
Since $(M,g,f)$ is a a steady gradient soliton, the sum $R+|\nabla f|^2$ is constant. Since $R \geq 0$ by Proposition \ref{nonnegative.scalar.curvature}, we conclude that $|\nabla f|$ is uniformly bounded from above. This completes the proof of Corollary \ref{bound.for.gradient.of.f}. \\

In the remainder of this section, we assume that $M$ is three-dimensional. As in Section \ref{background}, we denote by $\lambda_1 \leq \lambda_2 \leq \lambda_3$ the eigenvalues of the tensor $R \, g_{ij} - 2 \, \text{\rm Ric}_{ij}$. Then $R = \lambda_1+\lambda_2+\lambda_3$. Since $R \geq 0$, it follows that $\lambda_3 \geq 0$ at each point on $M$. 

\begin{proposition}[B.L.~Chen \cite{Chen}]
\label{lower.bound.for.lambda_1}
Assume that $n=3$. Then $k \lambda_1 + 2R \geq 0$ for every nonnegative integer $k$.
\end{proposition} 

\textbf{Proof.} 
The proof is by induction on $k$. For $k=0$, the assertion follows from Proposition \ref{nonnegative.scalar.curvature}. 

We now turn to the inductive step. Suppose that $k \geq 1$ and $(k-1)\lambda_1 + 2R \geq 0$. As above, we fix an arbitrary point $p \in M$. Let $N_0$ and $N_1$ denote the constants in Lemma \ref{barrier}. Let us fix a radius $r > 0$. We define a continuous function $v: B(p,r) \to \mathbb{R}$ by 
\[v(x) := k\lambda_1(x)+2R(x) + 4k(12+N_0+N_1 r) r^2 (r^2-d(p,x)^2)^{-2}\] 
for $x \in B(p,r)$. We claim that $v(x) \geq 0$ for all $x \in B(p,r)$. To prove this, we argue by contradiction. Let $\bar{x}$ be a point in $B(p,r)$ where the function $v$ attains its minimum, and suppose that $v(\bar{x}) < 0$. The evolution equation for the scalar curvature implies 
\[\Delta R + \langle X,\nabla R \rangle = -2 \, |\text{\rm Ric}|^2.\] 
The evolution equation for the Ricci tensor gives 
\[\Delta \lambda_1 + \langle X,\nabla \lambda_1 \rangle \leq -(\lambda_1^2 + \lambda_2 \lambda_3),\] 
where the inequality is understood in the barrier sense. More precisely, we can find a smooth function $\psi$ which is defined in an open neighborhood of $\bar{x}$ and satisfies the following conditions: 
\begin{itemize} 
\item $\psi(x) \geq \lambda_1(x)$ in an open neighborhood of the point $\bar{x}$.
\item $\psi(\bar{x}) = \lambda_1(\bar{x})$.
\item $\Delta \psi + \langle X,\nabla \psi \rangle \leq -(\lambda_1^2 + \lambda_2 \lambda_3)$ at the point $\bar{x}$. 
\end{itemize} 
For example, we may define $\psi := \text{\rm Ric}(\xi,\xi)$, where $\xi$ is a smooth unit vector field satisfying $\text{\rm Ric}(\xi,\xi) = \lambda_1$ at $\bar{x}$; $D\xi = 0$ at $\bar{x}$; and $\Delta \xi = 0$ at $\bar{x}$. 

By Lemma \ref{barrier}, we can find a smooth function $\rho$ which is defined in an open neighborhood of $\bar{x}$ and satisfies the following conditions: 
\begin{itemize}
\item $\rho(x) \geq d(p,x)^2$ in an open neighborhood of the point $\bar{x}$.
\item $\rho(\bar{x}) = d(p,\bar{x})^2$.
\item $|\nabla \rho|^2 = 4\rho$ at the point $\bar{x}$.
\item $\Delta \rho + \langle X,\nabla \rho \rangle \leq N_0 + N_1 \sqrt{\rho}$ at the point $\bar{x}$.
\end{itemize}
Then 
\[k\psi(x)+2R(x) + 4k(12+N_0+N_1 r) r^2 (r^2-\rho(x))^{-2} \geq v(x) \geq v(\bar{x})\] 
in an open neighborhood of $\bar{x}$, with equality at the point $\bar{x}$. Consequently, the function $k\psi+2R + 4k(12+N_0+N_1 r) r^2 (r^2-\rho)^{-2}$ attains a local minimum at the point $\bar{x}$. Thus, we conclude that 
\begin{align*} 
0 &\leq k \, (\Delta \psi + \langle X,\nabla \psi \rangle) + 2 \, (\Delta R + \langle X,\nabla R \rangle) \\ 
&+ 4k(12+N_0+N_1 r) r^2 \, \big [ \Delta((r^2-\rho)^{-2}) + \langle X,\nabla((r^2-\rho)^{-2}) \rangle \big ] \\ 
&= k \, (\Delta \psi + \langle X,\nabla \psi \rangle) + 2 \, (\Delta R + \langle X,\nabla R \rangle) \\ 
&+ 8k(12+N_0+N_1 r) r^2 (r^2-\rho)^{-3} \, (\Delta \rho + \langle X,\nabla \rho \rangle) \\ 
&+ 24k(12+N_0+N_1 r) r^2 (r^2-\rho)^{-4} \, |\nabla \rho|^2 \\ 
&\leq -k \, (\lambda_1^2+\lambda_2\lambda_3) - 4 \, |\text{\rm Ric}|^2 \\ 
&+ 8k(12+N_0+N_1 r) (N_0+N_1 \sqrt{\rho}) r^2 (r^2-\rho)^{-3} \\ 
&+ 96k(12+N_0+N_1 r) \rho r^2 (r^2-\rho)^{-4} 
\end{align*} 
at the point $\bar{x}$. Since $v \leq 0$ at the point $\bar{x}$, we know that 
\[-(k\lambda_1+2R) \geq 4k(12+N_0+N_1 r) r^2 (r^2-\rho)^{-2}\]  
at the point $\bar{x}$. Note that $R \geq 0$, $(k-1)\lambda_1 + 2R \geq 0$, $k\lambda_1 + 2R \leq 0$, $\lambda_1 \leq 0$, and $\lambda_3 \geq 0$ at the point $\bar{x}$. This implies 
\begin{align*} 
&k^2 \, (\lambda_1^2+\lambda_2\lambda_3) + 4k \, |\text{\rm Ric}|^2 \\ 
&= k^2 \, (\lambda_1^2+\lambda_2\lambda_3) + 2k \, (\lambda_1^2+\lambda_2^2+\lambda_3^2+\lambda_1\lambda_2+\lambda_2\lambda_3+\lambda_3\lambda_1) \\ 
&= (k\lambda_1+2R)^2 - 2R(k\lambda_1+2R) + k ((k-1)\lambda_1+2R) \lambda_3 \\ 
&+ k^2 (\lambda_2-\lambda_1)\lambda_3 - k\lambda_1\lambda_3 + 2k\lambda_2^2 \\ 
&\geq (k\lambda_1+2R)^2 \\  
&\geq 16k^2(12+N_0+N_1 r)^2 r^4 (r^2-\rho)^{-4} 
\end{align*} 
at the point $\bar{x}$. Putting these facts together, we obtain 
\begin{align*} 
0 &\leq -k \, (\lambda_1^2+\lambda_2\lambda_3) - 4 \, |\text{\rm Ric}|^2 \\ 
&+ 8k(12+N_0+N_1 r) (N_0+N_1 \sqrt{\rho}) r^2 (r^2-\rho)^{-3} \\ 
&+ 96k(12+N_0+N_1 r) \rho r^2 (r^2-\rho)^{-4} \\ 
&\leq -16k(12+N_0+N_1 r)^2 r^4 (r^2-\rho)^{-4} \\ 
&+ 8k(12+N_0+N_1 r) (N_0+N_1 r) r^4 (r^2-\rho)^{-4} \\ 
&+ 96k(12+N_0+N_1 r) r^4 (r^2-\rho)^{-4} \\ 
&= -8k(12+N_0+N_1 r)^2 r^4 (r^2-\rho)^{-4}
\end{align*} 
at the point $\bar{x}$. This is a contradiction. 

Thus, we conclude that 
\[k\lambda_1(x)+2R(x) + 4k(12+N_0+N_1 r) r^2 (r^2-d(p,x)^2)^{-2} \geq 0\] 
for all $x \in B(p,r)$. Sending $r \to \infty$ gives $k\lambda_1(x)+2R(x) \geq 0$ for each point $x \in M$. This completes the proof of Proposition \ref{lower.bound.for.lambda_1}. \\

\begin{corollary}[B.L.~Chen \cite{Chen}]
\label{nonnegative.curvature}
Assume that $n=3$. Then $(M,g)$ has nonnegative sectional curvature.
\end{corollary}

\textbf{Proof.} 
Sending $k \to \infty$ in Proposition \ref{lower.bound.for.lambda_1} gives $\lambda_1 \geq 0$. \\

\begin{corollary}
\label{bounded.curvature}
Assume that $n=3$. Then $(M,g)$ has bounded curvature. 
\end{corollary}

\textbf{Proof.}
Since $(M,g,f)$ is a a steady gradient soliton, the sum $R+|\nabla f|^2$ is constant. Consequently, the scalar curvature is uniformly bounded from above. Hence, the assertion follows from Corollary \ref{nonnegative.curvature}. \\

From now on, we will assume that $(M,g)$ is $\kappa$-noncollapsed. Moreover, we will assume that $(M,g,f)$ is normalized so that $R+|\nabla f|^2 = 1$ at each point on $M$. Let $\Phi_t$ denote the one-parameter group of diffeomorphisms generated by the vector field $-X$. It follows from Corollary \ref{bound.for.gradient.of.f} that $\Phi_t$ is defined for all $t \in (-\infty,\infty)$. In view of Corollary \ref{nonnegative.curvature} and Corollary \ref{bounded.curvature}, the metrics $\Phi_t^*(g)$, $t \in (-\infty,0]$, form an ancient $\kappa$-solution to the Ricci flow. \\

\begin{proposition}
\label{positive.curvature}
Assume that $n=3$ and $(M,g)$ is $\kappa$-noncollapsed. Then $(M,g)$ has positive sectional curvature.
\end{proposition}

\textbf{Proof.}
By Corollary \ref{nonnegative.curvature}, $(M,g)$ has nonnegative sectional curvature. We claim that $(M,g)$ has strictly positive sectional curvature. Suppose this is false. By the strict maximum principle, the universal cover of $(M,g)$ splits off a line. Using Perelman's classification of ancient $\kappa$-solutions in dimension $2$ (see Theorem \ref{2D.noncollapsed.case}), we conclude that the universal cover of $(M,g)$ is isometric to a cylinder $S^2 \times \mathbb{R}$, up to scaling. In particular, $(M,g)$ has constant scalar curvature. This contradicts the fact that $\Delta R + \langle X,\nabla R \rangle = -2 \, |\text{\rm Ric}|^2$ at each point on $M$. This completes the proof of Proposition \ref{positive.curvature}. \\

\begin{proposition}
\label{limit.of.rescalings.of.a.soliton}
Assume that $n=3$ and $(M,g)$ is $\kappa$-noncollapsed. Let $p_j$ be a sequence of points going to infinity, and let $r_j^{-2} := R(p_j)$. Let us dilate the manifold $(M,g)$ around the point $p_j$ by the factor $r_j^{-1}$. Then, after passing to a subsequence if necessary, the rescaled manifolds converge in the Cheeger-Gromov sense to a cylinder of radius $\sqrt{2}$. 
\end{proposition}

\textbf{Proof.} 
Since $(M,g)$ has positive sectional curvature, the assertion follows from Corollary \ref{limit.of.rescalings.of.an.ancient.solution}. \\

\begin{corollary}
\label{R.tends.to.0}
Assume that $n=3$ and $(M,g)$ is $\kappa$-noncollapsed. Then $R \to 0$ at infinity.
\end{corollary}

\textbf{Proof.} 
Proposition \ref{limit.of.rescalings.of.a.soliton} implies that $R^{-2} \, \Delta R \to 0$ and $R^{-\frac{3}{2}} \, |\nabla R| \to 0$ at infinity. Since $R$ is bounded from above, it follows that $\Delta R \to 0$ and $|\nabla R| \to 0$ at infinity. Since $|X|$ is bounded, we conclude that $\Delta R + \langle X,\nabla R \rangle \to 0$ at infinity. On the other hand, the evolution equation for the scalar curvature gives $\Delta R + \langle X,\nabla R \rangle = -2 \, |\text{\rm Ric}|^2$ at each point in $M$. Putting these facts together, we conclude that $|\text{\rm Ric}|^2 \to 0$ at infinity. This completes the proof of Corollary \ref{R.tends.to.0}. \\

\begin{corollary}
\label{limit.of.rescalings.of.a.soliton.2}
Assume that $n=3$ and $(M,g)$ is $\kappa$-noncollapsed. Let $p_j$ be a sequence of points going to infinity, and let $r_j^{-2} := R(p_j)$. Let us dilate the manifold $(M,g)$ around the point $p_j$ by the factor $r_j^{-1}$. Then, after passing to a subsequence if necessary, the rescaled manifolds converge in the Cheeger-Gromov sense to a cylinder of radius $\sqrt{2}$, and the rescaled vector fields $r_j \, X$ converge in $C_{\text{\rm loc}}^\infty$ to the axial vector field on the cylinder.
\end{corollary}

\textbf{Proof.} 
The vector field $X$ satisfies the pointwise estimates $|X| \leq 1$ and $|DX| = |\text{\rm Ric}| \leq C \, R$. Moreover, Perelman's pointwise derivative estimate (see Theorem \ref{pointwise.derivative.estimate}) implies $|D^{m+1} X| = |D^m \text{\rm Ric}| \leq C \, R^{\frac{m+2}{2}}$ for every positive integer $m$. Consequently, the rescaled vector fields $r_j \, X$ converge in $C_{\text{\rm loc}}^\infty$ to a limit vector field on the cylinder. Since $|X|^2 = 1-R \to 1$ at infinity, the limiting vector field on the cylinder has unit length at each point. Since $|DX| \leq C \, R$, the limiting vector field on the cylinder is parallel. This completes the proof of Corollary \ref{limit.of.rescalings.of.a.soliton.2}. \\

\begin{proposition} 
\label{critical.points.of.f}
Assume that $n=3$ and $(M,g)$ is $\kappa$-noncollapsed. Then the function $f$ has a unique critical point $p_*$, and $f$ attains its global minimum at the point $p_*$. 
\end{proposition}

\textbf{Proof.} 
In view of Corollary \ref{R.tends.to.0}, there exists a point $p_*$ where the scalar curvature is maximal. In particular, $\nabla R = 0$ at the point $p_*$. Since $R+|X|^2$ is constant, we know that $\nabla R + 2 \, \text{\rm Ric}(X) = 0$ at each point on $M$. Consequently, $\text{\rm Ric}(X) = 0$ at the point $p_*$. Since $(M,g)$ has positive Ricci curvature, it follows that $X=0$ at the point $p_*$. In other words, $p_*$ is a critical point of $f$. Since $(M,g)$ has positive Ricci curvature, the function $f$ is strictly convex. Thus, $p_*$ is the only critical point of $f$, and $f$ attains its global minimum at the point $p_*$. This completes the proof of Proposition \ref{critical.points.of.f}. \\

\begin{corollary} 
\label{linear.growth.of.f}
Assume that $n=3$ and $(M,g)$ is $\kappa$-noncollapsed. Then there exists a positive constant $C$ such that 
\[\frac{1}{C} \, d(p_*,x) \leq f(x) \leq C \, d(p_*,x)\] 
outside some compact set. 
\end{corollary}

\textbf{Proof.} 
The upper bound for $f$ follows from Corollary \ref{bound.for.gradient.of.f}. The lower bound follows from Proposition \ref{critical.points.of.f} together with the strict convexity of $f$. \\

\begin{proposition}[H.~Guo \cite{Guo}]
\label{asymptotic.behavior.of.scalar.curvature}
Assume that $n=3$ and $(M,g)$ is $\kappa$-noncollapsed. Then $f \, R \to 1$ at infinity.
\end{proposition}

\textbf{Proof.} 
Using the evolution equation for the scalar curvature, we obtain 
\[\Delta R + \langle X,\nabla R \rangle = -2 \, |\text{\rm Ric}|^2\] 
at each point in $M$. This implies 
\begin{align*} 
\langle X,\nabla (R^{-1}-f) \rangle 
&= -R^{-2} \, \langle X,\nabla R \rangle - |\nabla f|^2 \\ 
&= R^{-2} \, \Delta R + 2R^{-2} \, |\text{\rm Ric}|^2 - |\nabla f|^2 
\end{align*}
at each point in $M$. Using Proposition \ref{limit.of.rescalings.of.a.soliton}, we obtain $R^{-2} \, \Delta R \to 0$ and $R^{-2} \, |\text{\rm Ric}|^2 \to \frac{1}{2}$ at infinity. Moreover, Corollary \ref{R.tends.to.0} implies $|\nabla f|^2 = 1-R \to 1$ at infinity. Putting these facts together, we conclude that 
\[\langle X,\nabla (R^{-1}-f) \rangle \to 0\] 
at infinity. Hence, if $\varepsilon>0$ is given, then 
\[\langle X,\nabla (R^{-1}-(1+\varepsilon)f) \rangle \leq 0\] 
and 
\[\langle X,\nabla (R^{-1}-(1-\varepsilon)f) \rangle \geq 0\] 
outside a compact set. Integrating these inequalities along the integral curves of $X$, we obtain 
\[\sup_M (R^{-1}-(1+\varepsilon)f) < \infty\] 
and 
\[\inf_M (R^{-1}-(1-\varepsilon)f) > -\infty.\] 
Since $\varepsilon>0$ is arbitrary, we conclude that $f \, R \to 1$ at infinity. This completes the proof of Proposition \ref{asymptotic.behavior.of.scalar.curvature}. \\

In particular, if $s$ is sufficiently large and $f(x)=s$, then $x$ lies at center of a neck, and the radius of the neck is $(1+o(1)) \sqrt{2s}$.

\section{Rotational symmetry of noncollapsed steady gradient Ricci solitons in dimension $3$ -- the proof of Theorem \ref{soliton.classification}}

\label{symmetry.of.solitons}

Throughout this section, we assume that $(M,g,f)$ is a non-flat steady gradient Ricci soliton in dimension $3$ which is $\kappa$-noncollapsed. Moreover, we assume that $(M,g,f)$ is normalized so that $R+|\nabla f|^2 = 1$ at each point on $M$. It follows from Proposition \ref{critical.points.of.f} that $f$ is bounded from below. After adding a constant to $f$ if necessary, we may assume that $f$ is positive at each point on $M$. For abbreviation, we put $X := \nabla f$.  Let $\Phi_t$ denote the one-parameter group of diffeomorphisms generated by the vector field $-X$. By Corollary \ref{bound.for.gradient.of.f}, $\Phi_t$ is defined for all $t \in (-\infty,\infty)$. Moreover, the metrics $\Phi_t^*(g)$, $t \in (-\infty,0]$, form an ancient $\kappa$-solution to the Ricci flow. 

Let us fix a large real number $L$ and a small positive real number $\varepsilon_1$ so that the conclusion of the Neck Improvement Theorem holds. In view of Proposition \ref{asymptotic.behavior.of.scalar.curvature}, we can find a large constant $\Lambda$ with the following properties: 
\begin{itemize} 
\item $f(x) \, R(x,0)^{\frac{1}{2}} \geq 10^6 L$ for each point $x \in M$ with $f(x) \geq \frac{\Lambda}{2}$. 
\item If $f(\bar{x}) \geq \frac{\Lambda}{2}$, then $(\bar{x},0)$ lies at the center of an evolving $\varepsilon_1^2$-neck.
\end{itemize}
By a repeated application of the Neck Improvement Theorem, we obtain the following result.

\begin{proposition}
\label{iteration}
Suppose that $j$ is a nonnegative integer and $x$ is a point in $M$ with $f(x) \geq 2^{\frac{j}{400}} \, \Lambda$. Then the point $(x,0)$ is $2^{-j} \, \varepsilon_1$-symmetric.
\end{proposition}

\textbf{Proof.} 
The proof is by induction on $j$. For $j=0$, the assertion is true by our choice of $\Lambda$. Suppose next that $j \geq 1$, and the assertion is true for $j-1$. Let us fix a point $\bar{x} \in M$ satisfying $f(\bar{x}) \geq 2^{\frac{j}{400}} \, \Lambda$, and let $r^{-2} := R(\bar{x},0)$. By our choice of $\Lambda$, $f(\bar{x}) \, r^{-1} = f(\bar{x}) \, R(\bar{x},0)^{\frac{1}{2}} \geq 10^6 L$. Since $|\nabla f| \leq 1$ at each point on $M$, we obtain 
\[f(x) \geq f(\bar{x}) - Lr \geq (1-10^{-6}) \, f(\bar{x}) \geq (1-10^{-6}) \, 2^{\frac{j}{400}} \, \Lambda \geq 2^{\frac{j-1}{400}} \, \Lambda\] 
for each point $x \in B_g(\bar{x},Lr)$. This implies 
\[f(\Phi_t(x)) \geq f(x) \geq 2^{\frac{j-1}{400}} \, \Lambda\] for each point $x \in B_g(\bar{x},Lr)$ and each $t \leq 0$. We now apply the induction hypothesis. Hence, if $x \in B_g(\bar{x},Lr)$ and $t \leq 0$, then the point $(\Phi_t(x),0)$ is $2^{-j+1} \, \varepsilon_1$-symmetric. Since we are working on a self-similar solution, the point $(x,t)$ plays the same role as the point $(\Phi_t(x),0)$. Consequently, if $x \in B_g(\bar{x},Lr)$ and $t \leq 0$, then the point $(x,t)$ is $2^{-j+1} \, \varepsilon_1$-symmetric. Using the Neck Improvement Theorem, we conclude that the point $(\bar{x},0)$ is $2^{-j} \, \varepsilon_1$-symmetric. This completes the proof of Proposition \ref{iteration}. \\

\begin{corollary}
\label{almost.killing.vector.fields}
If $m$ is sufficiently large, then we can find smooth vector fields $U^{(1,m)},U^{(2,m)},U^{(3,m)}$ on the domain $\{m - 80 \sqrt{m} \leq f \leq m + 80 \sqrt{m}\}$ with the following properties: 
\begin{itemize}
\item $\sum_{l=0}^2 \sum_{a=1}^3 m^l \, |D^l(\mathscr{L}_{U^{(a,m)}}(g))|^2 \leq C \, m^{-800}$.
\item If $\Sigma \subset \{m - 80 \sqrt{m} \leq f \leq m + 80 \sqrt{m}\}$ is a leaf of the CMC foliation, then $\sup_\Sigma \sum_{a=1}^3 m^{-1} \, |\langle U^{(a)},\nu \rangle|^2 \leq C \, m^{-800}$, where $\nu$ denotes the unit normal vector to $\Sigma$.
\item If $\Sigma \subset \{m - 80 \sqrt{m} \leq f \leq m + 80 \sqrt{m}\}$ is a leaf of the CMC foliation, then 
\[\sum_{a,b=1}^3 \bigg | \delta_{ab} - \text{\rm area}(\Sigma)^{-2} \int_\Sigma \langle U^{(a,m)},U^{(b,m)} \rangle \, d\mu \bigg |^2 \leq C \, m^{-800}.\]
\end{itemize}
\end{corollary}

As explained in Section 7 of \cite{Brendle4}, we can glue approximate Killing vector fields on overlapping necks. This allows us to draw the following conclusion:

\begin{corollary}
\label{almost.killing.vector.fields.2}
We can find smooth vector fields $U^{(1)},U^{(2)},U^{(3)}$ such that 
\[|\mathscr{L}_{U^{(a)}}(g)| \leq C \, (f+100)^{-100},\] 
\[|D(\mathscr{L}_{U^{(a)}}(g))| \leq C \, (f+100)^{-100},\] 
\[|D^2(\mathscr{L}_{U^{(a)}}(g))| \leq C \, (f+100)^{-100}.\] 
Finally, given any positive real number $\varepsilon$, we have  
\[\sum_{a,b=1}^3 \bigg | \delta_{ab} - \text{\rm area}(\Sigma)^{-2} \int_\Sigma \langle U^{(a)},U^{(b)} \rangle \, d\mu \bigg |^2 \leq \varepsilon^2\]
whenever $\Sigma$ is a leaf of the CMC foliation which is sufficiently far out near infinity (depending on $\varepsilon$).
\end{corollary}

\begin{lemma} 
\label{U.almost.orthogonal.to.X}
We have $|\langle [U^{(a)},X],X \rangle| \leq C \, (f+100)^{-40}$.
\end{lemma}

\textbf{Proof.} 
The vector field $X$ satisfies $|X|^2 = 1-R$. Let us take the Lie derivative along $U^{(a)}$ on both sides. This gives 
\[(\mathscr{L}_{U^{(a)}}(g))(X,X) + 2 \, \langle \mathscr{L}_{U^{(a)}}(X),X \rangle = -\mathscr{L}_{U^{(a)}}(R).\] 
Using the formula for the linearization of the scalar curvature (see \cite{Besse}, Theorem 1.174 (e)), we obtain 
\[|\mathscr{L}_{U^{(a)}}(R)| \leq C \, |D^2(\mathscr{L}_{U^{(a)}}(g))| + C \, |\text{\rm Ric}| \, |\mathscr{L}_{U^{(a)}}(g)| \leq C \, (f+100)^{-40}.\] 
Putting these facts together, we conclude that 
\[|\langle \mathscr{L}_{U^{(a)}}(X),X \rangle| \leq C \, (f+100)^{-40}.\] 
This completes the proof of Lemma \ref{U.almost.orthogonal.to.X}. \\

\begin{lemma}
\label{[U,X].almost.parallel}
We have $|D([U^{(a)},X])| \leq C \, (f+100)^{-40}$.
\end{lemma}

\textbf{Proof.} 
The vector field $X$ satisfies $g_{jk} \, D_i X^j = \text{\rm Ric}_{ik}$. Let us take the Lie derivative along $U^{(a)}$ on both sides. Using the formula for the linearization of the Levi-Civita connection (see \cite{Besse}, Theorem 1.174 (a)), we obtain  
\begin{align*} 
&(\mathscr{L}_{U^{(a)}}(g))_{jk} \, D_i X^j + g_{jk} \, D_i(\mathscr{L}_{U^{(a)}}(X))^j \\ 
&+ \frac{1}{2} \, D_i (\mathscr{L}_{U^{(a)}}(g))_{jk} \, X^j + \frac{1}{2} \, D_j (\mathscr{L}_{U^{(a)}}(g))_{ik} \, X^j - \frac{1}{2} \, D_k (\mathscr{L}_{U^{(a)}}(g))_{ij} \, X^j \\ 
&= (\mathscr{L}_{U^{(a)}}(\text{\rm Ric}))_{ik}.
\end{align*} 
The formula for the linearization of the Ricci tensor (see \cite{Besse}, Theorem 1.174 (d)) gives 
\[|\mathscr{L}_{U^{(a)}}(\text{\rm Ric})| \leq C \, |D^2(\mathscr{L}_{U^{(a)}}(g))| + C \, |\text{\rm Rm}| \, |\mathscr{L}_{U^{(a)}}(g)| \leq C \, (f+100)^{-40}.\] 
Putting these facts together, we conclude that 
\[|D(\mathscr{L}_{U^{(a)}}(X))| \leq C \, (f+100)^{-40}.\] 
This completes the proof of Lemma \ref{[U,X].almost.parallel}. \\

\begin{lemma}
\label{estimate.for.Lie.bracket}
We have $|[U^{(a)},X]| \leq (f+100)^{-20}$ outside a compact set.
\end{lemma}

\textbf{Proof.} 
Suppose that the assertion is false. Then there exists an index $a \in \{1,2,3\}$ and a sequence of points $p_j$ going to infinity such that 
\[|[U^{(a)},X]| \geq (f+100)^{-20}\] 
at the point $p_j$. Let us define $s_j := f(p_j)$, and let $A_j$ denote the norm of the vector field $[U^{(a)},X]$ at the point $p_j$. By assumption, $A_j \geq (s_j+100)^{-20}$ for each $j$. Since $|\nabla f| \leq 1$ at each point on $M$, we know that $\frac{s_j}{2} \leq f \leq \frac{3s_j}{2}$ at each point in $B_g(p_j,\frac{s_j}{2})$. Using Lemma \ref{U.almost.orthogonal.to.X} and Lemma \ref{[U,X].almost.parallel}, we obtain 
\[\sup_{B_g(p_j,\frac{s_j}{2})} |\langle [U^{(a)},X],X \rangle| \leq C \, (s_j+100)^{-40} \leq C \, (s_j+100)^{-20} \, A_j\] 
and 
\[\sup_{B_g(p_j,\frac{s_j}{2})} |D([U^{(a)},X])| \leq C \, (s_j+100)^{-40} \leq C \, (s_j+100)^{-20} \, A_j.\] 
In the next step, we integrate the bound for $D([U^{(a)},X])$ along geodesics emanating from $p_j$. If $j$ is sufficiently large, we obtain 
\[\sup_{B_g(p_j,\frac{s_j}{2})} |[U^{(a)},X]| \leq A_j + C \, (s_j+100)^{-30} \leq 2A_j.\] 
We now dilate the manifold $(M,g)$ around the point $p_j$ by the factor $s_j^{-\frac{1}{2}}$. By Corollary \ref{limit.of.rescalings.of.a.soliton.2}, the rescaled manifolds converge in the Cheeger-Gromov sense to a cylinder of radius $\sqrt{2}$, and the rescaled vector fields $s_j^{\frac{1}{2}} \, X$ converge in $C_{\text{\rm loc}}^\infty$ to the axial vector field on the cylinder. Moreover, the vector fields $s_j^{\frac{1}{2}} \, A_j^{-1} \, [U^{(a)},X]$ converge in $C_{\text{\rm loc}}^{\frac{1}{2}}$ to a non-trivial parallel vector field on the cylinder, and this limiting vector field is orthogonal to the axial vector field on the cylinder. This is a contradiction. This completes the proof of Lemma \ref{estimate.for.Lie.bracket}. \\

\begin{lemma} 
\label{error.term.in.PDE}
We have $|\Delta U^{(a)} + D_X U^{(a)}| \leq C \, (f+100)^{-20}$.
\end{lemma}

\textbf{Proof.} 
Using Proposition \ref{divergence.of.Lie.derivative} and Corollary \ref{almost.killing.vector.fields.2}, we obtain 
\[|\Delta U^{(a)} + \text{\rm Ric}(U^{(a)})| \leq C \, |D(\mathscr{L}_{U^{(a)}}(g))| \leq C \, (f+100)^{-100}.\] 
Moreover, Lemma \ref{estimate.for.Lie.bracket} gives 
\[|[U^{(a)},X]| \leq C \, (f+100)^{-20}.\] 
Using the identity 
\begin{align*} 
\Delta U^{(a)} + D_X U^{(a)} 
&= \Delta U^{(a)} + D_{U^{(a)}} X - [U^{(a)},X] \\ 
&= \Delta U^{(a)} + \text{\rm Ric}(U^{(a)}) - [U^{(a)},X], 
\end{align*} 
we conclude that 
\[|\Delta U^{(a)} + D_X U^{(a)}| \leq C \, (f+100)^{-20}.\] 
This completes the proof of Lemma \ref{error.term.in.PDE}. \\

For abbreviation, we define smooth vector fields $Q^{(1)},Q^{(2)},Q^{(3)}$ by $Q^{(a)} := \Delta U^{(a)} + D_X U^{(a)}$.

\begin{proposition}
\label{existence.of.W}
We can find smooth vector fields $W^{(1)},W^{(2)},W^{(3)}$ such that $|W^{(a)}| \leq C \, (f+100)^{-8}$ and $\Delta W^{(a)} + D_X W^{(a)} = Q^{(a)}$. 
\end{proposition}

\textbf{Proof.}
We consider a sequence of real numbers $s_j \to \infty$. For each $j$ and each $a \in \{1,2,3\}$, we denote by $W^{(a,j)}$ the solution of the elliptic PDE 
\[\Delta W^{(a,j)} + D_X W^{(a,j)} = Q^{(a)}\] 
on the domain $\{f \leq s_j\}$ with Dirichlet boundary condition $W^{(a,j)} = 0$ on the boundary $\{f=s_j\}$. This Dirichlet problem has a solution by the Fredholm alternative. Moreover, since $Q^{(a)}$ is smooth, it follows that $W^{(a,j)}$ is smooth.

By Lemma \ref{error.term.in.PDE}, $Q^{(a)}$ satisfies a pointwise estimate of the form 
\[|Q^{(a)}| \leq K \, (f+100)^{-20},\] 
where $K$ is a large constant that does not depend on $j$. Using Kato's inequality, we obtain 
\[\Delta |W^{(a,j)}| + \langle X,\nabla |W^{(a,j)}| \rangle \geq -|Q^{(a)}| \geq -K \, (f+100)^{-20}\] 
on the set $\{W^{(a,j)} \neq 0\}$. On the other hand, using the identity $\Delta f + \langle X,\nabla f \rangle = R+|\nabla f|^2 = 1$, we obtain 
\begin{align*} 
&\Delta((f+100)^{-8}) + \langle X,\nabla ((f+100)^{-8}) \rangle \\ 
&= -8 \, (f+100)^{-9} \, (\Delta f + \langle X,\nabla f \rangle) + 72 \, (f+100)^{-10} \, |\nabla f|^2 \\ 
&\leq -8 \, (f+100)^{-9} + 72 \, (f+100)^{-10} \\ 
&\leq -(f+100)^{-9}. 
\end{align*} 
Using the maximum principle, we conclude that 
\[|W^{(a,j)}| \leq K \, (f+100)^{-8}\] 
on the set $\{f \leq s_j\}$. 

We now send $j \to \infty$. After passing to a subsequence, the vector fields $W^{(a,j)}$ converge in $C_{\text{\rm loc}}^\infty$ to a smooth vector field $W^{(a)}$. The limiting vector field $W^{(a)}$ satisfies 
\[|W^{(a)}| \leq K \, (f+100)^{-8}\] 
and 
\[\Delta W^{(a)} + D_X W^{(a)} = Q^{(a)}.\] 
This completes the proof of Proposition \ref{existence.of.W}. \\

\begin{proposition}
\label{bound.for.DW}
We have $|DW^{(a)}| \leq C \, (f+100)^{-8}$.
\end{proposition}

\textbf{Proof.} 
By Proposition \ref{existence.of.W}, the vector field $W^{(a)}$ satisfies $\Delta W^{(a)} + D_X W^{(a)} = Q^{(a)}$. This equation can be rewritten as $\Delta W^{(a)} + \mathscr{L}_X(W^{(a)}) + \text{\rm Ric}(W^{(a)}) = Q^{(a)}$. We next consider the vector fields $\Phi_t^*(W^{(a)})$ on the evolving background $(M,\Phi_t^*(g))$. These vector fields satisfy the parabolic PDE 
\[\frac{\partial}{\partial t} \Phi_t^*(W^{(a)}) = \Delta_{\Phi_t^*(g)} \Phi_t^*(W^{(a)}) + \text{\rm Ric}_{\Phi_t^*(g)}(\Phi_t^*(W^{(a)})) - \Phi_t^*(Q^{(a)}).\] 
The assertion follows now from interior estimates for parabolic PDE (see e.g. \cite{Brendle-Naff}, Proposition C.2). This completes the proof of Proposition \ref{bound.for.DW}. \\

We next define smooth vector fields $V^{(1)},V^{(2)},V^{(3)}$ by $V^{(a)} := U^{(a)} - W^{(a)}$.

\begin{proposition}
\label{pde.for.V}
The vector field $V^{(a)}$ satisfies $\Delta V^{(a)} + D_X V^{(a)} = 0$.
\end{proposition}

\textbf{Proof.} 
This follows immediately from Proposition \ref{existence.of.W}. \\

\begin{proposition}
\label{V.is.an.exact.Killing.vector.field}
The tensor $\mathscr{L}_{V^{(a)}}(g)$ vanishes identically.
\end{proposition} 

\textbf{Proof.} 
Recall that 
\[\Delta V^{(a)} + D_X V^{(a)} = 0.\] 
By Corollary \ref{evolution.of.Lie.derivative.soliton.version}, the tensor $h^{(a)} := \mathscr{L}_{V^{(a)}}(g)$ satisfies 
\[\Delta_L h^{(a)} + \mathscr{L}_X(h^{(a)}) = 0.\] 
Hence, Corollary \ref{anderson.chow.estimate.soliton.version} implies 
\[\Delta \Big ( \frac{|h^{(a)}|^2}{R^2} \Big ) + \Big \langle X,\nabla \Big ( \frac{|h^{(a)}|^2}{R^2} \Big ) \Big \rangle + \frac{2}{R} \, \Big \langle \nabla R,\nabla \Big ( \frac{|h^{(a)}|^2}{R^2} \Big ) \Big \rangle \geq 0.\] 
Using the maximum principle, we obtain 
\[\sup_{\{f \leq s\}} \frac{|h^{(a)}|}{R} \leq \sup_{\{f=s\}} \frac{|h^{(a)}|}{R}\] 
for each $s$. On the other hand, Corollary \ref{almost.killing.vector.fields.2} and Proposition \ref{bound.for.DW} imply 
\[|h^{(a)}| \leq |\mathscr{L}_{U^{(a)}}(g)| + C \, |DW^{(a)}| \leq C \, (f+100)^{-8}.\] 
Using Proposition \ref{asymptotic.behavior.of.scalar.curvature}, we deduce that 
\[\frac{|h^{(a)}|}{R} \leq C \, (f+100)^{-7}.\] 
In particular, 
\[\sup_{\{f=s\}} \frac{|h^{(a)}|}{R} \to 0\] 
as $s \to \infty$. Putting these facts together, we conclude that $h^{(a)}$ vanishes identically. This completes the proof of Proposition \ref{V.is.an.exact.Killing.vector.field}. \\

\begin{corollary}
\label{V.orthogonal.to.X}
We have $[V^{(a)},X] = 0$ and $\langle V^{(a)},X \rangle = 0$.
\end{corollary}

\textbf{Proof.} 
Note that 
\[\Delta V^{(a)} + D_X V^{(a)} = 0\] 
by Proposition \ref{pde.for.V}. On the other hand, since $V^{(a)}$ is a Killing vector field, we obtain 
\[\Delta V^{(a)} + \text{\rm Ric}(V^{(a)}) = 0\] 
by Proposition \ref{divergence.of.Lie.derivative}. This implies $D_X V^{(a)} = \text{\rm Ric}(V^{(a)})$. On the other hand, $D_{V^{(a)}} X = \text{\rm Ric}(V^{(a)})$. Consequently, $[V^{(a)},X] = 0$. This proves the first statement. We now turn to the proof of the second statement. Since $V^{(a)}$ is a Killing vector field, we obtain 
\[\nabla(\mathscr{L}_{V^{(a)}}(f)) = \mathscr{L}_{V^{(a)}}(\nabla f) = \mathscr{L}_{V^{(a)}}(X) = 0.\] 
Consequently, the function $\mathscr{L}_{V^{(a)}}(f) = \langle V^{(a)},X \rangle$ is constant. On the other hand, by Proposition \ref{critical.points.of.f}, the vector field $X$ vanishes at some point $p_* \in M$. Thus, we conclude that the function $\langle V^{(a)},X \rangle$ vanishes identically. This completes the proof of Corollary \ref{V.orthogonal.to.X}. \\

Corollary \ref{V.orthogonal.to.X} implies that the vector fields $V^{(1)},V^{(2)},V^{(3)}$ are tangential to the level sets of $f$.

\begin{proposition}
\label{orthonormality}
Given any positive real number $\varepsilon$, we have 
\[\sum_{a,b=1}^3 \bigg | \delta_{ab} - \text{\rm area}(\Sigma)^{-2} \int_\Sigma \langle V^{(a)},V^{(b)} \rangle \, d\mu \bigg |^2 \leq \varepsilon^2\]
whenever $\Sigma$ is a leaf of the CMC foliation which is sufficiently far out near infinity (depending on $\varepsilon$).
\end{proposition}

\textbf{Proof.} 
This follows by combining Corollary \ref{almost.killing.vector.fields.2} with Proposition \ref{existence.of.W}. \\

Proposition \ref{orthonormality} ensures that the vector fields $V^{(1)},V^{(2)},V^{(3)}$ are non-trivial near infinity. From this, it is easy to see that $(M,g)$ is rotationally symmetric.

\end{document}